\documentclass[12pt]{article}
\usepackage{xspace}
\usepackage[english,french]{babel}
\usepackage[latin1]{inputenc}
\usepackage[T1]{fontenc}
\usepackage{pslatex}
\usepackage{amsmath,amssymb,amsfonts,amsthm}
\usepackage{lscape}
\usepackage{url}
\usepackage{array}
\usepackage{graphicx}
\usepackage{hyperref}
\usepackage{listings}
\usepackage{color}
\usepackage{tikz-cd}
\usepackage{fullpage}
\usepackage{caption}

\addto\captionsfrench{}
\def\ZZ{\ensuremath{\mathbb Z}}
\def\QQ{\ensuremath{\mathbb Q}}
\def\RR{\ensuremath{\mathbb R}}
\def\PP{\ensuremath{\mathbb P}}

\def\cA{\ensuremath{\mathcal A}}
\def\cF{\ensuremath{\mathcal F}}
\def\cH{\ensuremath{\mathcal H}}
\def\cP{\ensuremath{\mathcal P}}
\def\cU{\ensuremath{\mathcal U}}
\def\cV{\ensuremath{\mathcal V}}
\def\cW{\ensuremath{\mathcal W}}

\def\smallmat#1{\left(\begin{smallmatrix}#1\end{smallmatrix}\right)}
\def\kbd#1{\texttt{#1}}
\DeclareMathOperator{\trace}{trace}
\DeclareMathOperator{\Sl}{SL}
\DeclareMathOperator{\Gl}{GL}
\DeclareMathOperator{\Psl}{PSL}

\DeclareMathOperator{\Hom}{Hom}

\theoremstyle{plain}
\newtheorem{thm}{Théorème}[section]
\newtheorem{lem}[thm]{Lemme}

\newtheorem{cor}[thm]{Corollaire}

\theoremstyle{definition}
\newtheorem{defn}[thm]{Définition}
\newtheorem{rem}[thm]{Remarque}
\newtheorem*{const}{Construction}
\newtheorem*{ex}{Exemple}

\def\farey{symbole de Farey étendu\xspace}
\def\muell{\mu_{ell}}
\def\aaa{\alpha}

\catcode`\@=11
\def\upparenfill{$\m@th\bracelu\leaders\vrule\hfill\braceru$}
\def\underparen#1{\mathop{\vtop{\ialign{##\crcr
      $\hfil\displaystyle{#1}\hfil$\crcr
      \noalign{\kern3\p@\nointerlineskip}
      \upparenfill\crcr\noalign{\kern3\p@}}}}\limits}

\def\encadre#1{{\color{red}{\underline{#1}}}}
\def\encadrebis#1{{\color{blue}{\underline{#1}}}}
\usepackage[Algorithme]{algorithm}
\usepackage[noend]{algorithmic}

\usepackage{multido}

\newcommand {\geo}[2] {(#2:1) arc (90+#2:270+#1:{tan((#2-#1)/2)})}
\newcommand {\sgeo}[2]{(#2:1) -- (#1:1)}
\newcommand{\unarc}[3]{({#1 * #3}:1.3) node {#2}}
\newcommand{\cell}[2]{({#1 * #2}:0.95) circle(0.05)}
\newcommand {\link}[5] {\unarc{#1}{#2}{#5} \geo {#1*#5}{#3*#5} \unarc{#3}{#4}{#5}}
\newcommand {\slink}[5] {\unarc{#1}{#2}{#5} \sgeo {#1*#5}{#3*#5}\unarc{#3}{#4}{#5}}

\title{Polygones fondamentaux d'une courbe modulaire}

\usepackage[noblocks]{authblk}

\author[*]{Karim Belabas}
\author[**]{Dominique Bernardi}
\author[***]{Bernadette Perrin-Riou}

\affil[*]{Univ. Bordeaux, CNRS, INRIA, IMB, UMR 5251, F-33400 Talence, France}
\affil[**]{Sorbonne Université, Institut de Mathématiques de Jussieu - Paris Rive Gauche, F-75005 Paris, France}
\affil[***]{Laboratoire de Mathématiques d'Orsay, Univ. Paris-Sud, CNRS, Université Paris-Saclay,
F-91405 Orsay, France.}

\begin{document}
\maketitle
\selectlanguage{english}
\begin{abstract}
A few pages in Siegel \cite{siegel} p. 115 (\S 2) describe how, starting with a fundamental
polygon for a compact Riemann surface, one can construct a symplectic basis of its homology.
This note retells that construction, specializing to the case where the surface is associated to
a subgroup $\Gamma$ of finite index in $\Psl_2(\ZZ)$. One then obtains by classical procedures
a generating system for $\Gamma$ with a minimal number of hyperbolic elements
and a presentation of the $\ZZ[\Gamma]$-module of the elements of $\ZZ[\PP^1(\QQ)]$
of degree 0.
\end{abstract}
\selectlanguage{french}
\begin{abstract}
Quelques pages de Siegel \cite{siegel} p. 115 (\S 2)
décrivent la construction d'une base symplectique de l'homologie d'une
surface de Riemann compacte à partir d'un polygone fondamental.
Cette note reprend cette construction
en l'appliquant au cas de la surface de Riemann associée
à un sous-groupe d'indice fini $\Gamma$ de $\Psl_2(\ZZ)$.
On en déduit par des procédés classiques un système de générateurs indépendants
de $\Gamma$ ayant un nombre minimal d'éléments hyperboliques et
une présentation du $\ZZ[\Gamma]$-module des éléments de $\ZZ[\PP^1(\QQ)]$ de degré 0.
\end{abstract}

%%On appelle domaine fondamental dans le
%%demi-plan de Poincaré complété $\cH \cup \PP^1(\QQ)$ un polygone hyperbolique
%%$P$ simplement connexe dont le bord est formé (essentiellement) de chemins géodésiques entre
%%des points de $\PP^1(\QQ)$ tel que 1) si $z$ est dans l'intérieur de $P$ et
%%$\gamma\in \Gamma$ est tel que $\gamma z\in P$, alors $\gamma = \Id$; 2) pour
%%chaque $z\in \cH$, il existe $\gamma\in\Gamma$ tel que $\gamma z \in P$.

Soit $\cH^*=\cH \cup \PP^1(\QQ)$ le demi-plan de Poincaré complété
et soit $\Gamma$ un sous-groupe d'indice fini de $\Psl_2(\ZZ)$ tel que le test
d'appartenance $\gamma \in \Gamma$ soit calculable, par exemple un
sous-groupe de congruence.
Il lui est associé (de manière non unique) un symbole de Farey qui par
définition est essentiellement
la donnée d'un polygone hyperbolique convexe de $\cH^*$ dont les sommets sont dans $\PP^1(\QQ)$
et d'une involution sur les chemins géodésiques de son bord
(voir la définition \ref{def:farey} pour un énoncé précis), ce qui permet de trouver
un domaine fondamental dans $\cH^*$ pour le quotient $X(\Gamma)=\Gamma\backslash \cH^*$.

De tels symboles fournissent des outils pour manipuler
algorithmiquement le groupe $\Gamma$, par exemple un système de
générateurs indépendants ou un test déterminant si $\Gamma$ est un sous-groupe
de congruence ainsi que son niveau. Ils admettent une construction effective:
les algorithmes principaux sont dus à Kulkarni~\cite{kulkarni} développant
des idées remontant à Poincaré~\cite{poincare}.
%%et reposent sur la construction de symboles de Farey associés au groupe
%% $\Gamma$.

%%Tout sous-groupe d'indice fini $\Gamma$ de
%%$\Psl_2(\ZZ)$ peut être défini par un symbole de Farey qui détermine un
%%domaine fondamental pour $\Gamma$ et un système de générateurs indépendants.

Les invariants usuels comme le nombre de points elliptiques, le nombre de
pointes ou le genre $g$ de la courbe modulaire $X(\Gamma)$ se lisent
directement sur le symbole. Il suffit d'un test effectif d'appartenance à
$\Gamma$ pour le construire, mais il fournit aussi réciproquement un tel
test. Ces algorithmes sont disponibles dans plusieurs systèmes de calcul
formel. Sans prétendre être exhaustif, citons une implantation dans Sage
basée sur le paquet \kbd{KFarey}~\cite{KurtLong08}, le paquet
\kbd{Congruence}~\cite{congruence} pour Gap et une implantation dans Magma
adaptée de Verrill~\cite{verrill}. Pour $\Gamma$ un sous-groupe de
congruence, la motivation principale est la manipulation des espaces de
formes modulaires associées en calculant des espaces de symboles modulaires
et l'action de l'algèbre de Hecke,
voir par exemple Cremona~\cite{cremona}. Le lien entre symboles modulaires et
symboles de Farey est repris de façon éclairante par Pollack et
Stevens~\cite{PS}, qui ramènent le calcul d'un espace de symboles modulaires
$\Hom_\Gamma(\ZZ[\PP^1(\QQ)]_0, V)$ pour une représentation $V$ de $\Gamma$
au calcul de la structure du groupe $\ZZ[\PP^1(\QQ)]_0$ des diviseurs
de degré $0$ sur $\PP^1(\QQ)$ comme $\Gamma$-module, laquelle est indépendante
de $V$ et est une conséquence immédiate de l'algorithme de Kulkarni.

Dans cette note, nous nous éloignons des motivations usuelles des
espaces de symboles modulaires ou de la représentation efficace des sous-groupes
de congruences et définissons une notion de symbole de Farey \emph{normalisé}
qui fournit un système de générateurs indépendants de $\Gamma$ contenant
exactement $2g$ éléments hyperboliques et une base symplectique pour
$H_1(X(\Gamma), \QQ)$. Nous donnons un algorithme pour passer d'un symbole de
Farey à un symbole de Farey normalisé, ce qui fournit une dissection (ou
représentation) canonique de la surface de Riemann $X(\Gamma)$. Précisons
d'emblée que cette notion de symbole normalisé n'a à notre connaissance pas
d'intérêt algorithmique : l'algorithme de normalisation fournit des
générateurs de taille bien plus importante que ceux associés au symbole
utilisé en entrée, sans permettre d'opération nouvelle.

Dans la première partie, nous faisons des rappels sur les symboles de Farey et
introduisons la notion de symbole de Farey normalisé. Nous rappelons le lien avec
les polygones fondamentaux de la surface de Riemann $X(\Gamma)$ associée à un
sous-groupe $\Gamma$ de $\Psl_2(\ZZ)$ ainsi qu'une application connue sur
la structure de $\ZZ[\PP^1(\QQ)]_0$ en tant que $\Gamma$-module, qui reste
valable en utilisant les symboles de Farey étendus et qui est fondamentale
dans la construction de symboles modulaires de Pollack et Stevens.

Dans la seconde partie, nous donnons l'algorithme de calcul d'un symbole de
Farey normalisé à partir d'un symbole de Farey quelconque. Dans l'appendice,
nous donnons des exemples calculés dans Pari/GP, système dans lequel nous
avons fait l'implantation de cet algorithme dans le cas particulier
du sous-groupe de Hecke $\Gamma_0(N)$.

Nous remercions le rapporteur pour ses remarques qui nous ont permis
de rendre plus clair le contenu de ce texte (nous l'espérons en tout cas).

\tableofcontents

\section{Polygone fondamental associé à un sous-groupe de \texorpdfstring{$\Psl_2(\ZZ)$}{PSL2(Z)}}
\label{farey}
% !TEX root = ../fareydissection.tex
% !TEX encoding = IsoLatin
On note $\cH$ le demi-plan de Poincaré et $\cH^*=\cH \sqcup \PP^1(\QQ)$
le demi-plan complété.
\subsection{Lemme préliminaire}
Si $r$ et $s$ sont deux éléments distincts de $\PP^1(\QQ)$, on note
$a=(r,s)$ l'arc géodésique dans $\cH^*$ menant de $r$ à $s$, et
$A_a=\smallmat{x_a&y_a\\z_a&t_a}$
une matrice de $M_2(\ZZ)$ primitive et de déterminant strictement positif
représentant~$a$ : on a donc $r=\frac{x_a}{z_a}$ et $s=\frac{y_a}{t_a}$
avec $\gcd(x_a,z_a)= \gcd(y_a,t_a)=1$.
La matrice $A_a$ est unique au signe près.
Son déterminant est appelé la largeur de l'arc $a$.
La matrice $$A_{a}^{-}=\begin{pmatrix}y_a&-x_a\\t_a&-z_a\end{pmatrix}
=A_a \begin{pmatrix}0&-1\\ 1&0\end{pmatrix}$$
représente l'arc $\overline{a}=(s,r)$ opposé à $a$.

%ordre sur $\PP^1(\QQ)$, Ou bien on prend l'ordre circulaire ???
%On conviendra que $\infty$ est plus grand que tout élément de $\QQ$ dans $\PP^1(\QQ)$.
On dit que des points $r_1, \dots, r_n$ de $\PP^1(\QQ)$ avec $n\geq 3$ sont dans l'ordre
circulaire si leurs images $R_1,\dots, R_n$ dans le cercle unité par la
transformation de Cayley $z \mapsto \frac{i-z}{i+z}$ sont dans cet ordre
lorsqu'on parcourt le cercle dans le sens trigonométrique.
%%On note alors abusivement $r_1\leq \dots \leq r_n$.

\begin{lem}\label{lemme:hpe}
Soient $r$, $s$, $t$ et $u$ quatre éléments de $\PP^1(\QQ)$ tels que
$r\neq s$, $t\neq u$ et tels que $r,s,t,u$ soient dans l'ordre circulaire.
On pose $a=(r,s)$ et $a^*=(t,u)$.

\begin{enumerate}
\item La matrice $\gamma=A_{a}{(A_{a^*}^-)}^{-1}$ de $\Gl_2(\QQ)$
vérifie $a=\gamma \overline{a^*}$.
%%%$\gamma=A_{a^*}^- A_{a}^{-1} \in \Gl_2(\QQ)$
\begin{enumerate}
\item Si $s \neq t$, alors la matrice $\gamma$ est hyperbolique.
\item Si $s = t$, alors la matrice $\gamma$ est parabolique si et seulement
si $\gamma$ est de déterminant 1 et hyperbolique sinon.
\end{enumerate}
\item
La matrice $\gamma= A_a \smallmat{0&1\\-1&0}A_a^{-1}=
A_a (A_a^-)^{-1}$ est une matrice de
$\Sl_2(\QQ)$ vérifiant
$a = \gamma\overline{a}$. Elle est d'ordre 2 dans $\Psl_2(\QQ)$.
\item La matrice
%%%$$\gamma= A_a \begin{pmatrix}0&-1\\1&-1\end{pmatrix}A_a^{-1}=A_a^{-} \begin{pmatrix}1&-1\\0&1\end{pmatrix}A_a^{-1}$$
$$\gamma= A_a^-\begin{pmatrix}0&-1\\1&-1\end{pmatrix}(A_a^-)^{-1}=
A_a \begin{pmatrix}-1&1\\0&-1\end{pmatrix}(A_a^-)^{-1}$$
est d'ordre 3 et vérifie $1 + \gamma +\gamma^2 =0$. Les arcs
$(a,\gamma a ,\gamma^2 a)$ forment un chemin fermé.
\end{enumerate}
\end{lem}

\begin{proof}
Nous utiliserons les faits suivants :
si $\gamma$ appartient à $\Gl_2^+(\QQ)$, elle est
parabolique si $\trace(\gamma)^2=4\det(\gamma)$ et
hyperbolique si elle a deux points fixes dans $\PP^1(\RR)$.
Elle est elliptique si $\trace(\gamma)^2< 4\det(\gamma)$
c'est-à-dire si elle n'a pas de points fixes dans $\PP^1(\RR)$.

\begin{enumerate}
\item
On vérifie que dans chacun des cas, la matrice $\gamma$ vérifie la propriété
annoncée.
\begin{enumerate}
\item
Si $s \neq t$, montrons que $\gamma$ est hyperbolique. Elle induit une application
continue de $\PP^1(\RR)$ dans $\PP^1(\RR)$.
L'image par $\gamma^{-1}$ du segment $\lbrack r,s\rbrack$ de $\PP^1(\RR)$ est le segment
$[u,t]$ qui contient $[r,s]$ ; par le théorème
des valeurs intermédiaires, $\gamma^{-1}$ a un point fixe
dans le segment $\lbrack r,s\rbrack$, donc $\gamma$ a un point fixe dans
le segment $\lbrack r,s \rbrack$.
De même, l'image du segment $\lbrack t,u\rbrack$ par $\gamma$ est le segment
$[s,r]$ qui contient $\lbrack t,u\rbrack$, donc $\gamma$
a un point fixe dans le segment $[t,u]$. Ayant deux points fixes
dans $\PP^1(\RR)$, $\gamma$ est hyperbolique.

\item Si $s=t$, un calcul explicite montre que le déterminant de $\gamma$ est
$\det(A_{a})/\det(A_{a^*})$ et que sa trace est $1 + \det \gamma$. On en
déduit que $\trace(\gamma)^2 - 4\det(\gamma) = (1 - \det \gamma)^2$
et donc que $\gamma$ est parabolique si et seulement si elle est de
déterminant 1. Dans le cas contraire, comme $\gamma$ a un point fixe dans
$\PP^1(\QQ)$, elle est nécessairement hyperbolique.
\end{enumerate}
\item La vérification de (2) et (3) est un calcul explicite.
\end{enumerate}
\end{proof}
\begin{rem} Les formules dans le cas elliptique sont obtenues de la manière suivante.
\begin{enumerate}
\item
La matrice
$\smallmat{0&1\\-1&0}$ est d'ordre 2 et
transforme les points $\infty,i,0$ en $0,i,\infty$, respectivement. La matrice
$$\gamma=A_a \begin{pmatrix}0&1\\-1&0\end{pmatrix} A_a^{-1}
$$
transforme les points $r,A_a(i),s$ en $s,A_a(i),r$.
%Par un calcul explicite, $\gamma=(\det A_a)^{-1}\begin{pmatrix}*&*\\ z_a^2 + t_a^2&*\end{pmatrix}$.
\item
Posons $\rho=\frac{1}{2} + i\frac{\sqrt{3}}{2}$.
La matrice $\tau=\smallmat{0&-1\\1&-1}$ est d'ordre 3,
laisse fixe $\rho$ et envoie les points $\infty,0,1$ sur $0,1,\infty$.
Si $T$ est le triangle hyperbolique de sommets $(0,\rho,\infty)$,
le triangle hyperbolique de sommets $(\infty, 0,1)$ est fixe par $\tau$,
est réunion des trois triangles hyperboliques $T$, $\tau T$ et $\tau^2 T$, et
contient l'unique point fixe $\rho$ de $\tau$ dans $\cH^*$. Ainsi, la matrice
$\gamma=A_a^- \smallmat{0&-1\\1&-1} {(A_a^-)}^{-1}$ de $\Sl_2(\QQ)$ laisse
fixe le triangle hyperbolique $(s,r,A_a^-(1))$.
\end{enumerate}
\end{rem}
\subsection{Symbole de Farey étendu}
\begin{defn}\label{def:farey}
On appelle \textsl{\farey} la donnée $(\cV,*,\muell)$
\begin{enumerate}
\item[(F1)] d'un polygone hyperbolique convexe orienté dont les sommets
$(r_1,\dots,r_n)$ sont dans $\PP^1(\QQ)$; son bord
est donc une suite $\cV$ d'arcs consécutifs $(a_1, \dots, a_n)$
avec $a_1=(r_1,r_2)$, $a_2=(r_2,r_3)$, $\dots, a_n=(r_n,r_1)$; on suppose que
l'un des $a_i$ est l'arc $(\infty,0)$;
\item[(F2)] d'une involution $* \colon a\mapsto a^*$ de $\cV$;
\item[(F3)] d'une application $\muell$ de l'ensemble $\cV_{ell}$
des points fixes de $*$ à valeurs dans $\{2,3\}$;
\end{enumerate}
tels que
\begin{enumerate}
\item[(F4)] si $a \in \cV - \cV_{ell}$, la matrice
%%%$\gamma_a=A_{a^*}^- A^{-1}$
$\gamma_a=A_a (A_{a^*}^-)^{-1}$
est dans $\Sl_2(\ZZ)$ ;
\item[\hphantom{(F4)}] si $a \in \cV_{ell}$ et $\muell(a)=2$, la matrice
%%%$\gamma_a=A_a^- A_a^{-1}$
$\gamma_a=A_a (A_a^-)^{-1}$
est dans $\Sl_2(\ZZ)$ ;
\item[\hphantom{(F4)}] si $a \in \cV_{ell}$ et $\muell(a)=3$, la matrice
%%%$\gamma_a= A_a^{-} \begin{pmatrix}1&-1\\0&1\end{pmatrix}A_a^{-1}$
$\gamma_a= A_a^- \smallmat{0&-1\\1&1}(A_a^-)^{-1}$
est dans $\Sl_2(\ZZ)$.
\end{enumerate}
\end{defn}
On appelle les éléments de $\cV$ les \emph{arcs} du \farey et les $\gamma_{a}$
pour $a\in\cV$ ses \emph{données de recollement}. Combinatoirement,
nous considèrerons $\cV$ comme un \emph{collier}: c'est l'orbite du mot
$a_1\dots a_n$ sous l'action du groupe $\ZZ/n\ZZ$ des permutations
circulaires. Nous appellerons \emph{facteur} de $\cV$ un mot
$a_ia_{i+1}\dots a_{j-1}a_j$ formé de lettres « consécutives »,
en considérant $i,j$ dans $\ZZ/n\ZZ$.

Pour $n \in\{2,3\}$, on
note $\cV_{ell,n}$ l'image réciproque de $n$ par $\muell$. Si $a \in
\cV_{ell}$, on dit que $a$ est un arc elliptique d'\textsl{ordre}
$\muell(a)$. Ainsi,
\begin{itemize}
%%%$a^*= \gamma_{a} \overline{a}$
\item pour $a\notin \cV_{ell,3}$, on a $a=\gamma_{a} \overline{a^*}$;
\item pour $a\in \cV_{ell,3}$, le triangle hyperbolique
$(a,\gamma_a a,\gamma_a^2 a)$ est invariant par $\gamma_a$.
\end{itemize}
\begin{defn}
On appelle \textsl{groupe} du \farey le sous-groupe
$\Gamma\subset \Psl_2(\ZZ)$ engendré par les données de recollement $\gamma_a$, $a\in \cV$.
\end{defn}
\begin{defn}
Un \farey $\cF=(\cV,*,\muell)$ est \emph{unimodulaire} si les arcs
de $\cF$ sont de largeur 1, autrement dit si les
matrices $A_{a}$ pour $a\in \cV$ sont de déterminant 1.
La condition (F4) est alors trivialement vérifiée.
\end{defn}
Graphiquement, nous représenterons les arcs de $\cV$ par des points sur un
cercle et l'involution $*$ sur l'ensemble des arcs de $\cV$
en reliant par un trait un arc (représenté par un point) et son image par~$*$.
Les éléments de $\cV_{ell,3}$ (resp. $\cV_{ell,2}$) sont indiqués par un
point plein (resp. creux), voir figure \ref{fig:involution}.
\begin{figure}[H]
\begin{center}
\begin{tikzpicture}
\input fareydraw/cercle0\string_15.tex
\end{tikzpicture}
\quad
\begin{tikzpicture}
\input fareydraw/cercle0\string_37.tex
\end{tikzpicture}
\caption{Exemples de représentation de l'involution $*$}
\label{fig:involution}
\end{center}
\end{figure}
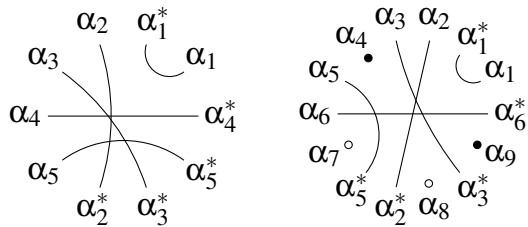
\begin{rem}
Ces notions sont introduites par Kulkarni~\cite{kulkarni}, qui appelle
simplement symbole de Farey ce que nous appelons symbole de Farey
unimodulaire. Dans la suite, nous
définirons un algorithme de normalisation qui abandonne l'unimodularité:
notre normalisation conservera la propriété $\gamma_a \in \Sl_2(\ZZ)$ pour
les données de recollement, sans qu'on ait nécessairement $A_a\in \Sl_2(\ZZ)$
pour chaque arc $a$ de $\cV$.
\end{rem}

Si $\cV=(a_1,\dots, a_n)$, on définit la distance de $a_i$ à $a_j$ par
$d(a_i,a_j) = \min(|j-i|, |n-i+j|)$.
Par exemple, $d(a_1,a_n)$ est égale à 1.
Le lemme suivant est une conséquence du lemme \ref{lemme:hpe}.
\begin{lem}
Soit $\cF=(\cV,*,\muell)$ un \farey de groupe $\Gamma$.
Si $d(a,a^*) \geq 2$, alors $\gamma_a$ est une matrice hyperbolique.
Si $d(a,a^*) = 1$, alors $\gamma_a$ est une matrice parabolique.
Si $d(a,a^*) = 0$, alors $\gamma_a$ est une matrice elliptique.
\end{lem}
On note $\cV_{hyp}$ (resp. $\cV_{par}$)
l'ensemble des arcs $a$ à distance $\geq 2$ (resp. 1) de leur image par $*$.
L'ensemble $\cV_{ell}$ déjà défini est l'ensemble des arcs à distance 0 de leur image
par $*$.

\begin{figure}[h]
\begin{center}
\begin{tikzpicture}[scale=15]
\fill [black!10!white] (21/100,1/8) -- (21/100,0) arc (180:0:1/50) arc (180:0:1/24) arc (180:0:1/30)
-- (2/5,0) -- (2/5,1/8);
\draw (21/100,0) arc (180:0:1/50) node (r) [below] {$r$} arc (180:100:1/24) node{$\circ$} node [below] {$t$} arc (100:0:1/24) node [below] (r) {$s$};
\draw (0.25,0) arc (180:100:1/24) node{$\circ$} arc (100:0:1/24);
\draw (1/3,0) arc (180:0:1/30);
\end{tikzpicture}
\begin{tikzpicture}[scale=15]
\fill [black!10!white] (1/4, 1/8) -- (1/4,0) arc (180:0:1/24) arc (180:73.1736:1/21) arc (133.1736:0:1/16) arc (180:0:1/16) -- (5/8,1/8) -- (1/4,1/8);
\draw (1/4,0) arc (180:0:1/24) node (r) {} node [below]{$r$};
\draw (r) arc (180:0:1/12) node [below]{$s$};
\draw (r) arc (180:73.1736:1/21) node (t) {$\bullet$} node [above] {$t$}
arc (133.1736:0:1/16) arc (0:180:1/12);
\draw (1/3,0) arc (180:0:1/30) node (z){} arc (180:0:1/20) arc (180:0:1/16);
\draw (2/5,0) arc (0:11.927:21/95);
\end{tikzpicture}
\end{center}
\caption{Points elliptiques sur le bord de polygones hyperboliques}
\label{fig:domaine}
\end{figure}
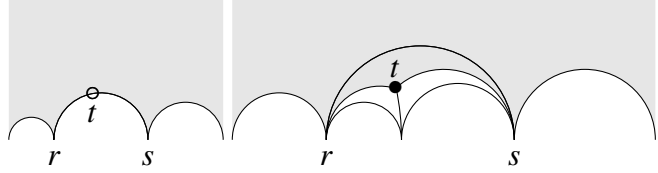

\begin{defn} Si $\cF=(\cV,*,\muell)$ est un \farey, on définit
le polygone hyperbolique $\cP(\cF)$ de $\cH^*$ dont les côtés sont
\begin{itemize}
\item les arcs géodésiques $a \in \cV$ tels que $a$ n'est pas un point fixe de $^*$ ;
\item les arcs géodésiques $(r,t)$ et $(t,s)$
où $(r,t,s)$ est l'image de $(0,i,\infty)$ (resp. de $(0,\rho,\infty)$)
par un élément de $\Psl_2(\QQ)$ et où $(r,s)$ appartient à $\cV_{ell,2}$
(resp. à $\cV_{ell,3}$).
\end{itemize}
On note $\mathcal{U(\cF)}$ l'enveloppe convexe de $\cP(\cF)$ dans $\cH^*$.
\end{defn}
Le polygone $\cP(\cF)$ est un polygone spécial au sens de \cite{kulkarni}
à condition d'enlever la condition d'unimodularité.
La figure \ref{fig:domaine} montre une partie d'un polygone hyperbolique
associé à un \farey ayant des arcs elliptiques et son enveloppe
convexe. Le symbole $\circ$ désigne un point elliptique d'ordre 2
et le symbole $\bullet$ désigne un point elliptique d'ordre 3.

\begin{rem}
Le choix de définir $\gamma_a$ par $a=\gamma_a \overline{a}^*$ (voir \cite{PS})
plutôt que par $\overline{a}^*=\gamma_a a$ (voir \cite{kulkarni}) permet de voir
$\gamma_a$ comme l'unique élément de $\Gamma$ modulo $\pm \text{Id}$ tel que
le côté $a$ du polygone $\cP(\cF)$ soit un côté commun de $\cP(\cF)$ et de $\gamma_a \cP(\cF)$.
\end{rem}
\begin{thm}\cite{kulkarni, maskit}
Soit $\Gamma$ un sous-groupe d'indice fini de $\Psl_2(\ZZ)$.
\begin{enumerate}
\item Il existe un symbole de Farey unimodulaire de groupe $\Gamma$.
\item Soit $\cF = (\cV, *,\muell)$ un \farey de groupe $\Gamma$ contenu dans
$\Psl_2(\ZZ)$.
Alors, $\cU(\cF)$ est un domaine fondamental admissible au sens de \cite{kulkarni}
pour l'action de $\Gamma$ sur $\cH$.
\end{enumerate}
\end{thm}
Les $\gamma_a$ pour $a\in \cV$ forment un système de générateurs
de $\Gamma$.
On dit aussi que $\cP(\cF)$ est un polygone fondamental pour $\Gamma$.

\begin{proof}
On peut trouver une démonstration de l'existence du symbole de Farey
unimodulaire $\cF$ de groupe $\Gamma$ dans \cite{kulkarni} (théorème 3.3).
%mais doit remonter à Poincaré (peut-être \cite{poincare}).
La deuxième partie se démontre comme le théorème 3.2 de \cite{kulkarni}.
\end{proof}
%%%Des implémentations ont été faites dans Magma, Sage et plus récemment dans Pari/GP (\cite{pari}).
%%%C'est cette dernière implémentation faite uniquement lorsque $\Gamma=\Gamma_0(N)$
%%%que nous utiliserons pour nos exemples (l'algorithme de construction.
%%%utilisé est celui que l'on trouve dans \cite{PS}).

\begin{rem}
Si $s$ est l'extrémité d'un arc $a = (s,t)$ de $ \cV - \cV_{ell}$ et si $\gamma_a$ est
la donnée de recollement associée à cet arc $a$ d'origine $s$, on appelle
successeur de $s$ le point $\gamma_a^{-1} s$.
Soit \hbox{$(s_0=s,\cdots, s_{l-1})$} l'orbite de $s$ :
$s_{j+1}=\gamma_{a_{j}}^{-1}s_j$, avec $a_{j}=(s_j,t_j)\in \cV$, et
$s_{0}=\gamma_{a_{l-1}}^{-1}s_{l-1}$.
Alors, le produit $\gamma_{a_{l-1}}^{-1}\cdots \gamma_{a_0}^{-1}$
est un générateur du stabilisateur de $s$ dans
$\Gamma$ et est conjugué à $\smallmat{1&w\\0&1}$
où $w>0$ est la largeur de la pointe~$s$.
En effet, posons
$a=a_0=(s,t)$; la matrice $\delta=\gamma_{a_{l-1}}^{-1}\cdots \gamma_{a_0}^{-1}$
envoie $s$ sur $s$.
Notons $r$ l'image de $t$ par $\delta$ et $b=(r,s) \in \cV$.
Nécessairement, $r$, $s$ et $t$ sont dans l'ordre circulaire
et $\delta$ vaut $A_b^- A_a^{-1}$ à un scalaire près.
Si $g$ est un élément de $\Sl_2(\QQ)$ qui envoie l'infini sur $s$,
par exemple $g=\smallmat{s&-1\\1&0}$, alors
$g^{-1}\delta g$ stabilise l'infini et vaut à un scalaire près
$$\begin{pmatrix}\det(A_a)&-y_r y_s + y_s y_t\\
0 &  \det(A_b)
\end{pmatrix},\quad\text{avec}\quad
A_a=\begin{pmatrix}x_s&x_t\\y_s&y_t \end{pmatrix}
\quad\text{et}\quad A_b=\begin{pmatrix}x_r&x_s\\y_r&y_s \end{pmatrix}.
$$
Comme $\det(A_a)=(s-t)y_sy_t>0$, on en déduit que $y_sy_t<0$ ;
de même $y_r y_s <0$. Ceci démontre que
$\delta=\gamma_{a_{l-1}}^{-1}\cdots \gamma_{a_0}^{-1}$ est conjugué
à
$\smallmat{1&w\\0&1}$ avec $w > 0$.
\end{rem}
\subsection{Symbole de Farey normalisé}

\begin{defn}Soit $\cF=(\cV,*, \muell)$ un \farey.
Deux éléments distincts $a$ et $b$ de $\cV$ sont
dits \emph{liés} si les quatre arcs $a$, $b$, $a^*$ et $b^*$
sont dans l'ordre $a$, $b$, $a^*$, $b^*$ ou $a$, $b^*$, $a^*$, $b$
dans $\cV$ (à permutation circulaire près).
Cette notion est symétrique et stable par $*$.
\end{defn}

\begin{defn}
Un \farey $\cF=(\cV,*,\muell)$ est dit \emph{désentrelacé} si
pour tout couple $(a, b)$ d'éléments liés,
les distances de $a$ à $a^*$ et de $b$ à $b^*$ sont égales à 2.
Il est dit \emph{normalisé} si pour tout élément $a$, la distance de $a$ à $a^*$
est inférieure ou égale à 2, l'égalité ayant lieu si et seulement si
$a$ est lié à un autre élément.
\end{defn}
Un \farey normalisé est désentrelacé.
\begin{ex}
\begin{enumerate}
\item Si $a$, $b$, $c$ et $d$ des arcs qui ne sont pas fixés par $*$,
$abca^*b^*c^*dd^*$ ou $acc^*ba^*b^*dd^*$ ne sont pas désentrelacés :
$a$ et $b$ sont liés mais $d(a,a^*)$ est égal à $3$ dans le premier cas et à $4$ dans
le second cas.
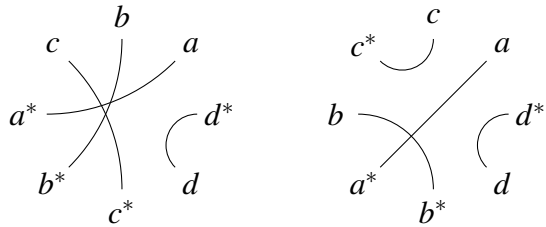
\begin{figure}[h]
\begin{center}
\begin{tikzpicture}[scale=1]
\draw \link {1}{$a$}{4}{$a^*$}{360/8};
\draw \link {2}{$b$}{5}{$b^*$}{360/8};
\draw \link {3}{$c$}{6}{$c^*$}{360/8};
\draw \link {7}{$d$}{8}{$d^*$}{360/8};
\end{tikzpicture}
\quad\quad
\begin{tikzpicture}[scale=1]
\draw \slink {1}{$a$}{5}{$a^*$}{360/8};
\draw \link {2}{$c$}{3}{$c^*$}{360/8};
\draw \link {4}{$b$}{6}{$b^*$}{360/8};
\draw \link {7}{$d$}{8}{$d^*$}{360/8};
\end{tikzpicture}
\end{center}
\caption{Involution de symboles non désentrelacés}
\end{figure}
\item Si $a$, $b$ et $c$ des arcs qui ne sont pas fixés par $*$
et $d$ un arc fixé par $*$,
$aba^*b^*cdc^*$ est désentrelacé et n'est pas normalisé car $d(c,c^*) =
2$, mais $c$ n'est pas lié à un autre élément.
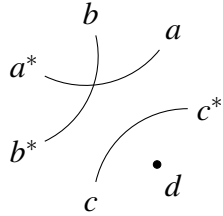
\begin{figure}[h]
\begin{center}\begin{tikzpicture}[scale=1]
\draw \link {1}{$a$}{3}{$a^*$}{360/7};
\draw \link {2}{$b$}{4}{$b^*$}{360/7};
\draw \link {5}{$c$}{7}{$c^*$}{360/7};
\draw \unarc{6}{$d$}{360/7};
\draw \cell{6}{360/7};
\fill \cell{6}{360/7};
\end{tikzpicture}
\end{center}
\caption{Involution d'un symbole désentrelacé mais non normalisé}
\end{figure}
\end{enumerate}
\end{ex}
Les exemples suivants se déduisent des calculs que nous détaillerons plus tard.
\begin{ex}[$\Gamma_0(15)$]
Le symbole de Farey unimodulaire
$$
1
\underset{ {\aaa_{1}}}{\underparen{}}\infty
\underset{ {\aaa^*_{1}}}{\underparen{}}0
\underset{ {\aaa_{2}}}{\underparen{}}\frac{1}{5}
\underset{ {\aaa_{3}}}{\underparen{}}\frac{1}{4}
\underset{ {\aaa_{4}}}{\underparen{}}\frac{1}{3}
\underset{ {\aaa_{5}}}{\underparen{}}\frac{2}{5}
\underset{ {\aaa^*_{2}}}{\underparen{}}\frac{1}{2}
\underset{ {\aaa^*_{3}}}{\underparen{}}\frac{3}{5}
\underset{ {\aaa^*_{5}}}{\underparen{}}\frac{2}{3}
\underset{ {\aaa^*_{4}}}{\underparen{}}1
$$
est associé à $\Gamma_0(15)$. Les arcs $\aaa_2$ et $\aaa_3$ sont liés
mais la distance entre $\aaa_2$ et $\aaa_2^*$ est $4 > 2$ donc
le \farey n'est pas désentrelacé. On a
$$\cV_{par}/*=\{\aaa_1\}, \cV_{hyp}/*=\{\aaa_2,\aaa_3,\aaa_4,\aaa_5\}, \cV_{ell}=\emptyset \ .$$
L'involution est représentée de la manière suivante:
\begin{center}
\begin{tikzpicture}
\input fareydraw/cercle0\string_15.tex
\end{tikzpicture}
\end{center}
Le \farey
$$
1
\underset{ {\aaa_{1}}}{\underparen{}}\infty
\underset{ {\aaa^*_{1}}}{\underparen{}}0
\underset{ {\aaa_{2}}}{\underparen{}}\frac{1}{12}
\underset{ {\aaa^*_{2}}}{\underparen{}}\frac{5}{59}
\underset{ {\aaa_{3}}}{\underparen{}}\frac{3}{35}
\underset{ {\aaa^*_{3}}}{\underparen{}}\frac{13}{151}
\underset{ {\aaa_{4}}}{\underparen{}}\frac{1}{8}
\underset{ {\aaa_{5}}}{\underparen{}}\frac{3}{13}
\underset{ {\aaa^*_{4}}}{\underparen{}}\frac{8}{19}
\underset{ {\aaa^*_{5}}}{\underparen{}}1

$$
est aussi associé à $\Gamma_0(15)$ et est normalisé.
Remarquons au passage que l'on voit sur ce \farey
qu'il y a 3 pointes non équivalentes autres que la pointe $0$,
qu'il n'y a pas de points elliptiques et que le genre de la courbe est 1 :
$$\cV^{norm}_{par}/*=\{\aaa_1,\aaa_2,\aaa_3\}, \cV^{norm}_{hyp}/*=\{\aaa_4,\aaa_5\},
\cV^{norm}_{ell}=\emptyset \ .$$

$${\color{red}{1}}
\underset{\aaa_1}{\underparen{}}\boxed{\infty}
\underset{\aaa_1^*}{\underparen{}}{\color{red}{0}}
\underset{\aaa_2}{\underparen{}}\boxed{\frac{1}{12}}
\underset{\aaa_2^*}{\underparen{}}{\color{red}{\frac{5}{59}}}
\underset{\aaa_3}{\underparen{}}\boxed{\frac{3}{35}}
\underset{\aaa_3^*}{\underparen{}}
\boxed{{\color{red}{\frac{13}{151}}}
\underset{\aaa_4}{\underparen{}}{\color{red}{\frac{1}{8}}}
\underset{\aaa_5}{\underparen{}}{\color{red}{\frac{3}{13}}}
\underset{\aaa_4^*}{\underparen{}}{\color{red}{\frac{8}{19}}}
\underset{\aaa_5^*}{\underparen{}}{\color{red}{1}}}
$$
\begin{center}
\begin{tikzpicture}
\input fareydraw/cercle1\string_15.tex
\end{tikzpicture}
\end{center}
\end{ex}
\begin{ex}[$\Gamma_0(20)$]
Le \farey suivant associé à $\Gamma_0(20)$ est unimodulaire mais non
normalisé :
$$
1
\underset{ {\aaa_{1}}}{\underparen{}}\infty
\underset{ {\aaa^*_{1}}}{\underparen{}}0
\underset{ {\aaa_{2}}}{\underparen{}}\frac{1}{5}
\underset{ {\aaa_{3}}}{\underparen{}}\frac{1}{4}
\underset{ {\aaa_{4}}}{\underparen{}}\frac{2}{7}
\underset{ {\aaa_{5}}}{\underparen{}}\frac{3}{10}
\underset{ {\aaa^*_{5}}}{\underparen{}}\frac{1}{3}
\underset{ {\aaa_{6}}}{\underparen{}}\frac{3}{8}
\underset{ {\aaa^*_{3}}}{\underparen{}}\frac{2}{5}
\underset{ {\aaa_{7}}}{\underparen{}}\frac{1}{2}
\underset{ {\aaa^*_{7}}}{\underparen{}}\frac{3}{5}
\underset{ {\aaa^*_{2}}}{\underparen{}}\frac{2}{3}
\underset{ {\aaa^*_{4}}}{\underparen{}}\frac{3}{4}
\underset{ {\aaa^*_{6}}}{\underparen{}}1

$$
\begin{center}
\begin{tikzpicture}
\draw \link {1}{$\aaa_{1}$}{2}{$\aaa^*_{1}$}{180/7};
\draw \link {3}{$\aaa_{2}$}{12}{$\aaa^*_{2}$}{180/7};
\draw \link {4}{$\aaa_{3}$}{9}{$\aaa^*_{3}$}{180/7};
\draw \link {5}{$\aaa_{4}$}{13}{$\aaa^*_{4}$}{180/7};
\draw \link {6}{$\aaa_{5}$}{7}{$\aaa^*_{5}$}{180/7};
\draw \link {8}{$\aaa_{6}$}{14}{$\aaa^*_{6}$}{180/7};
\draw \link {10}{$\aaa_{7}$}{11}{$\aaa^*_{7}$}{180/7};

\end{tikzpicture}
\end{center}
$$\cV_{par}/*=\{\aaa_1,\aaa_5, \aaa_7\},, \cV_{hyp}/*=\{\aaa_2, \aaa_3, \aaa_4,\aaa_6\}, \cV_{ell}=\emptyset \ .$$
Le \farey suivant associé à $\Gamma_0(20)$ est normalisé :
$$
1
\underset{ {\aaa_{1}}}{\underparen{}}\infty
\underset{ {\aaa^*_{1}}}{\underparen{}}0
\underset{ {\aaa_{2}}}{\underparen{}}\frac{1}{5}
\underset{ {\aaa^*_{2}}}{\underparen{}}\frac{4}{19}
\underset{ {\aaa_{3}}}{\underparen{}}\frac{3}{14}
\underset{ {\aaa^*_{3}}}{\underparen{}}\frac{11}{51}
\underset{ {\aaa_{4}}}{\underparen{}}\frac{107}{496}
\underset{ {\aaa^*_{4}}}{\underparen{}}\frac{524}{2429}
\underset{ {\aaa_{5}}}{\underparen{}}\frac{8}{37}
\underset{ {\aaa_{6}}}{\underparen{}}\frac{3}{11}
\underset{ {\aaa^*_{5}}}{\underparen{}}\frac{73}{253}
\underset{ {\aaa^*_{6}}}{\underparen{}}\frac{20}{69}
\underset{ {\aaa_{7}}}{\underparen{}}\frac{3}{10}
\underset{ {\aaa^*_{7}}}{\underparen{}}1

$$
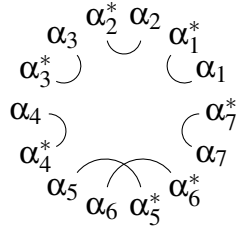
\begin{center}
\begin{tikzpicture}
\draw \link {1}{$\aaa_{1}$}{2}{$\aaa^*_{1}$}{180/7};
\draw \link {3}{$\aaa_{2}$}{4}{$\aaa^*_{2}$}{180/7};
\draw \link {5}{$\aaa_{3}$}{6}{$\aaa^*_{3}$}{180/7};
\draw \link {7}{$\aaa_{4}$}{8}{$\aaa^*_{4}$}{180/7};
\draw \link {9}{$\aaa_{5}$}{11}{$\aaa^*_{5}$}{180/7};
\draw \link {10}{$\aaa_{6}$}{12}{$\aaa^*_{6}$}{180/7};
\draw \link {13}{$\aaa_{7}$}{14}{$\aaa^*_{7}$}{180/7};

\end{tikzpicture}
\end{center}
$$\cV^{norm}_{par}/*=\{\aaa_1, \aaa_2, \aaa_3, \aaa_4, \aaa_7\}, \cV^{norm}_{hyp}/*=\{\aaa_5, \aaa_6\},
\cV^{norm}_{ell}=\emptyset \ .$$
On voit que la taille des rationnels intervenant dans le \farey a augmenté.
\end{ex}
\begin{ex}[$\Gamma_0(37)$]
Le symbole de Farey unimodulaire suivant pour $\Gamma_0(37)$
$$
1
\underset{ {\aaa_{1}}}{\underparen{}}\infty
\underset{ {\aaa^*_{1}}}{\underparen{}}0
\underset{ {\aaa_{2}}}{\underparen{}}\frac{1}{5}
\underset{ {\aaa_{3}}}{\underparen{}}\frac{1}{4}
\underset{\underset{\bullet}{\aaa_{4}}}{\underparen{}}\frac{1}{3}
\underset{ {\aaa_{5}}}{\underparen{}}\frac{3}{8}
\underset{ {\aaa_{6}}}{\underparen{}}\frac{2}{5}
\underset{\underset{\circ}{\aaa_{7}}}{\underparen{}}\frac{3}{7}
\underset{ {\aaa^*_{5}}}{\underparen{}}\frac{1}{2}
\underset{ {\aaa^*_{2}}}{\underparen{}}\frac{4}{7}
\underset{\underset{\circ}{\aaa_{8}}}{\underparen{}}\frac{3}{5}
\underset{ {\aaa^*_{3}}}{\underparen{}}\frac{2}{3}
\underset{\underset{\bullet}{\aaa_{9}}}{\underparen{}}\frac{3}{4}
\underset{ {\aaa^*_{6}}}{\underparen{}}1

$$
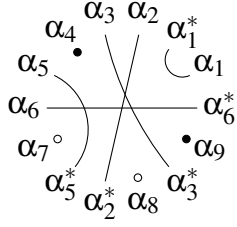
\begin{center}
\begin{tikzpicture}
\draw \link {1}{$\aaa_{1}$}{2}{$\aaa^*_{1}$}{180/7};
\draw \slink {3}{$\aaa_{2}$}{10}{$\aaa^*_{2}$}{180/7};
\draw \link {4}{$\aaa_{3}$}{12}{$\aaa^*_{3}$}{180/7};
\draw \unarc{5}{$\aaa_{4}$}{180/7};
\draw \cell{5}{180/7};
\fill \cell{5}{180/7};
\draw \link {6}{$\aaa_{5}$}{9}{$\aaa^*_{5}$}{180/7};
\draw \slink {7}{$\aaa_{6}$}{14}{$\aaa^*_{6}$}{180/7};
\draw \unarc{8}{$\aaa_{7}$}{180/7};
\draw \cell{8}{180/7};
\draw \unarc{11}{$\aaa_{8}$}{180/7};
\draw \cell{11}{180/7};
\draw \unarc{13}{$\aaa_{9}$}{180/7};
\draw \cell{13}{180/7};
\fill \cell{13}{180/7};
\end{tikzpicture}
\end{center}
représenté plus simplement par
$$\aaa_1 \ \aaa_1^*\ \aaa_2\ \aaa_3
\ \underset{\bullet}{\aaa_4}
\ \aaa_5 \ \aaa_6\
\underset{\circ}{\aaa_7}\ \aaa_5^*\ \aaa_2^*
\ \underset{\circ}{ \aaa_8} \ \aaa_3^*
\ \underset{\bullet}{\aaa_9}\ \aaa_6^*
$$
est un symbole de Farey associé à $\Gamma_0(37)$.
Il n'est pas désentrelacé car, par exemple, $\aaa_2$ et $\aaa_3$ sont liés mais la distance de $\aaa_2$
à $\aaa_2^*$ est strictement supérieure à 2.
Dans cet exemple, il n'y a que deux pointes non équivalentes.
Par contre, il n'y a que deux matrices paraboliques correspondant à $\aaa_1$
et à $\aaa_1^*$.
On a
$$\cV_{par}/*=\{\aaa_1\}, \cV_{hyp}/*=\{\aaa_2,\aaa_3,\aaa_5,\aaa_6\}, \cV_{ell}=\{\aaa_4,\aaa_7,\aaa_8,\aaa_9\}\ .$$
Un \farey normalisé obtenu par notre algorithme est
$$
1
\underset{ {\aaa_{1}}}{\underparen{}}\infty
\underset{ {\aaa^*_{1}}}{\underparen{}}0
\underset{\underset{\circ}{\aaa_{2}}}{\underparen{}}\frac{1}{6}
\underset{ {\aaa_{3}}}{\underparen{}}\frac{2}{11}
\underset{ {\aaa_{4}}}{\underparen{}}\frac{1}{5}
\underset{ {\aaa^*_{3}}}{\underparen{}}\frac{1}{4}
\underset{ {\aaa^*_{4}}}{\underparen{}}\frac{11}{43}
\underset{\underset{\bullet}{\aaa_{5}}}{\underparen{}}\frac{21}{82}
\underset{ {\aaa_{6}}}{\underparen{}}\frac{10}{39}
\underset{ {\aaa_{7}}}{\underparen{}}\frac{6}{23}
\underset{ {\aaa^*_{6}}}{\underparen{}}\frac{5}{19}
\underset{ {\aaa^*_{7}}}{\underparen{}}\frac{51}{193}
\underset{\underset{\circ}{\aaa_{8}}}{\underparen{}}\frac{7}{26}
\underset{\underset{\bullet}{\aaa_{9}}}{\underparen{}}1

$$
\begin{center}
\begin{tikzpicture}
\draw \link {1}{$\aaa_{1}$}{2}{$\aaa^*_{1}$}{180/7};
\draw \unarc{3}{$\aaa_{2}$}{180/7};
\draw \cell{3}{180/7};
\draw \link {4}{$\aaa_{3}$}{6}{$\aaa^*_{3}$}{180/7};
\draw \link {5}{$\aaa_{4}$}{7}{$\aaa^*_{4}$}{180/7};
\draw \unarc{8}{$\aaa_{5}$}{180/7};
\draw \cell{8}{180/7};
\fill \cell{8}{180/7};
\draw \link {9}{$\aaa_{6}$}{11}{$\aaa^*_{6}$}{180/7};
\draw \link {10}{$\aaa_{7}$}{12}{$\aaa^*_{7}$}{180/7};
\draw \unarc{13}{$\aaa_{8}$}{180/7};
\draw \cell{13}{180/7};
\draw \unarc{14}{$\aaa_{9}$}{180/7};
\draw \cell{14}{180/7};
\fill \cell{14}{180/7};
\end{tikzpicture}
\end{center}

On a
$$\cV^{norm}_{par}/*=\{\aaa_1\}, \cV^{norm}_{hyp}/*=\{\aaa_3,\aaa_4,\aaa_6,\aaa_7\},
\cV^{norm}_{ell}=\{\aaa_2,\aaa_5,\aaa_8,\aaa_9\}\ . $$
\end{ex}

\begin{lem}
Soit $\cF=(\cV,*,\muell)$ un \farey normalisé de groupe $\Gamma$.
\begin{enumerate}
\item Les extrémités des arcs qui sont à la distance 0 ou 2 de leur
image par $^*$ sont toutes équivalentes modulo $\Gamma$ ainsi
que les extrémités externes des facteurs $cc^*$ pour $c\in \cV_{par}$.
\item
Soit $P_0$ la classe d'équivalence commune.
L'application
$$\cV_{par}/* \to \Gamma\backslash \PP^1(\QQ) - \{P_0\}$$
qui associe à $c$ la classe de son extrémité commune avec $c^*$
est une bijection.
\end{enumerate}
\end{lem}
\begin{proof}
Le collier $\cV$ est formé de facteurs de la forme $aba^*b^*$, $cc^*$ et $d$ avec $d=d^*$.
Dans le cas $aba^*b^*$, on a
%%%$\gamma_a a^*=\overline{a}$ et $\gamma_b b^*=\overline{b}$.
$ \overline{a}=\gamma_a{a^*}$ et $ \overline{b}=\gamma_bb^*$.
Si $a=(r,s)$ et $b=(s,t)$, on a nécessairement
%%%$s= \gamma_a^{-1}\gamma_b^{-1}\gamma_a r$, $t=\gamma_b^{-1}\gamma_a r$ :
$$u= \gamma_a^{-1} r,\quad s= \gamma_a\gamma_b\gamma_a^{-1} r, \quad t=\gamma_a^{-1} \gamma_b r,
\quad
v=\gamma_b^{-1}\gamma_a\gamma_b \gamma_a^{-1} r\ . $$
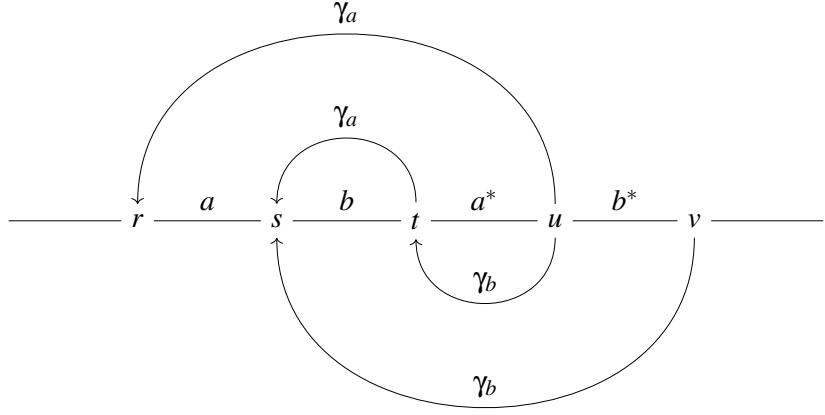
\begin{center}
\resizebox{0.7\textwidth}{!}{
\begin{tikzpicture}[node distance = 2cm]
\node (q) {};
\node (r) [right of= q] {$r$};
\node (s) [right of= r] {$s$};
\node (t) [right of= s] {$t$};
\node (u) [right of= t] {$u$};
\node (v) [right of= u] {$v$};
\node (w) [right of=v] {};
\path (q) edge (r);
\path (r) edge node [above] {$a$} (s);
\path (s) edge node [above] {$b$} (t);
\path (t) edge node [above] {$a^*$} (u);
\path (u) edge node [above] {$b^*$} (v);
\path (v) edge (w);

\draw[->] (u) .. node [above]{$\gamma_a$} controls +(up:3.5cm) and +(up:3.5cm) .. (r);
\draw[->] (v) .. node [above]{$\gamma_b$} controls +(down:3.5cm) and +(down:3.5cm) .. (s);
\draw[->] (t) .. node [above]{$\gamma_a$} controls +(up:1.5cm) and +(up:1.5cm) .. (s);
\draw[->] (u) .. node [above]{$\gamma_b$} controls +(down:1.5cm) and +(down:1.5cm) .. (t);
\end{tikzpicture}
}
\end{center}
Donc les extrémités des arcs $a$, $b$, $a^*$, $b^*$ sont toutes équivalentes modulo $\Gamma$.
%%%$$aba^*b^*=(r,\gamma_a^{-1}\gamma_b^{-1}\gamma_a r)\;
%%%(\gamma_a^{-1}\gamma_b^{-1}\gamma_a r,\gamma_b^{-1}\gamma_a r)\;
%%%(\gamma_b^{-1}\gamma_a r,\gamma_a r)\;
%%%(\gamma_a r,\gamma_b\gamma_a^{-1}\gamma_b^{-1}\gamma_a r)$$
Dans le cas d'un facteur de la forme $d$ avec $d=d^*$,
on a soit $d = \gamma_d \overline{d}$, soit
$d+\gamma_d d +\gamma_d^2 d=0$; dans les deux cas, les extrémités de $d$
sont équivalentes modulo $\Gamma$. Finalement, dans le cas d'un facteur de la
forme $cc^*$, les extrémités externes sont équivalentes puisque
%%%$\gamma_c c = \overline{c^*}$.
$\gamma_c^{-1} c = \overline{c^*}$.

On en déduit que toutes les extrémités externes de tous les facteurs de la
forme $aba^*b^*$, $cc^*$ et $d = d^*$ sont équivalentes ainsi que les
  extrémités internes des facteurs de la forme $aba^*b^*$. Cela termine la
  démonstration de la partie (1).

Il reste à regarder ce qui se passe pour la seconde extrémité de $c$ lorsque $c \in
\cV_{par}$ (facteur du type $cc^*$ avec $c\neq c^*$).
Les éléments de $\PP^1(\QQ)$ qui ne sont pas équivalents à $P_0$
sont nécessairement équivalents modulo $\Gamma$
à une des extrémités des arcs de $\cV$ et donc
à une extrémité de $c$ pour un $c$ de $\cV_{par}$.
Il suffit donc de démontrer que pour deux facteurs distincts du type $c_1c_1^*$ et $c_2c_2^*$,
 le point $s_1$ commun à $c_1$ et $c_1^*$ et le point
$s_2$ commun à $c_2$ et $c_2^*$ ne sont pas équivalents modulo $\Gamma$.
Prenons un petit disque ouvert dont le bord contient $s_1$ et inclus dans $\cH$.
Son image dans la surface $\Gamma\backslash \cH$ est un disque épointé
et ne peut pas contenir l'image d'un point proche de $s_2$.
Par continuité, $s_1$ et $s_2$ ne peuvent pas
être équivalents modulo $\Gamma$.
\end{proof}
\begin{rem}
Résumons les résultats précédents ainsi que des résultat connus
(voir par exemple \cite{kulkarni}).
Soit $g(\Gamma)$ le genre de $X(\Gamma)$, $\nu_\infty(\Gamma)$ le nombre de
pointes de $X(\Gamma)$, $\nu_2(\Gamma)$ et $\nu_3(\Gamma)$ le nombre
de points elliptiques d'ordre~2 et~3. Soit un \farey
 normalisé $\cF=(\cV,*,\muell)$. Alors, $\cV$ est formé de
 \begin{enumerate}
 \item $g(\Gamma)$ facteurs du type $a b{a^*}{b^*}$
 avec les $a$, $b$, ${a^*}$ et ${b^*}$ distincts deux à deux
 ($a$ et $b \in \cV_{hyp}$);
 \item $\nu_\infty(\Gamma)-1$ facteurs de type $c{c^*}$ avec $c\neq{c^*}$
 ($c\in \cV_{par}$) ;
 \item $\nu_2(\Gamma)+ \nu_3(\Gamma)$ points fixes $c$ par
 l'involution $*$ ($c \in \cV_{ell}$).
 \end{enumerate}
Le groupe fondamental $\pi_1(Y(\Gamma))$
  est libre de rang $r = 2g + t - 1$.
Le nombre d'arcs du domaine fondamental associé à $\cF$
est $2\left(r + \#\cV_{ell,2} + \#\cV_{ell,3}\right)$.
Le groupe $\Gamma$ est isomorphe à un produit libre de $\#\cV_{ell,2}$
exemplaires de $\ZZ/2\ZZ$, de $\#\cV_{ell,3}$ exemplaires de $\ZZ/3\ZZ$
et de $r$ exemplaires de $\ZZ$ (\cite{kulkarni}, prop. 3.4 et prop.7).
\end{rem}

\subsection{\texorpdfstring{Structure de $\ZZ[\PP^1(\QQ)]_0$}{Lg}}
\label{structure}
% !TEX root = ../fareydissection.tex
%!TEX encoding = IsoLatin

Soit $\Delta= \ZZ[\PP^1(\QQ)]$ le groupe des diviseurs sur $\PP^1(\QQ)$
et $\Delta_0$ le sous-groupe des diviseurs de degré 0.
L'action naturelle de $G = \Gl_2(\QQ)$ sur $\PP^1(\QQ)$ induit une action de $G$
sur $\Delta$.
Si $c \in \PP^1(\QQ)$, on note $\{c\}$ le diviseur associé à $c$ dans $\ZZ[\PP^1(\QQ)]$.
Si $c_1$ et $c_2$ sont dans $\PP^1(\QQ)$, on note $(c_1,c_2)$ le diviseur
$\{c_2\} - \{c_1\}$ associé dans $\Delta_0$. On parlera aussi du chemin $(c_1,c_2)$ :
dans l'interprétation géométrique, il s'agit du chemin géodésique
reliant les deux pointes $c_1$ et $c_2$ du demi-plan de Poincaré complété
$\cH^\ast$.
Le groupe $\Delta_0$ est engendré par les diviseurs de la forme
\hbox{$(c_1,c_2) = \{c_2\} - \{c_1\}$} avec $c_1$, $c_2 \in \PP^1(\QQ)$.
Les résultats de Pollack et Stevens se transposent de la manière suivante.
\begin{thm}
Soit $\cF$ un \farey $(\cV,*, \muell)$ et
$(\gamma_a)$ ses données de recollement pour $a\in \cV$.
Le $\ZZ[\Gamma]$-module $\Delta_0$ est défini par
générateurs et relations de la manière suivante :
l'image de $\cV$ dans $\Delta_0$ forme un système de générateurs vérifiant les relations
\begin{enumerate}
\item $\sum_{a\in \cV} a = 0$ ;
\item $a + \gamma_a a^*=0$ si $a \notin \cV_{ell,3}$ ;
\item $a+\gamma_a a +\gamma_a^2 a =0$ si $a \in \cV_{ell,3}$.
\end{enumerate}
\end{thm}
Pour la démonstration, voir \cite{PS}.
\def\EE{(\cV - \cV_{ell})/*}
On peut facilement extraire un système de générateurs minimal.
Soit $\EE$ un ensemble formé de représentants du quotient de $\cV - \cV_{ell}$
par l'involution~$*$. Le premier élément de $\EE$ est supposé être le chemin
$(\infty,0)$. On note $\cA = \EE \,\bigsqcup\, \cV_{ell}$.

\begin{cor}
Le $\ZZ[\Gamma]$-module $\Delta_0=\ZZ[\PP^1(\QQ)]_0$ est le quotient du
$\ZZ[\Gamma]$-module libre engendré par les éléments de $\cA$ par les relations
\begin{equation*}
\begin{cases}
\sum_{a\in \cA} \lambda_a a&=0 \\
\mu_a a &= 0
\text{\; pour $a \in \cV_{ell}$}
\end{cases}
\end{equation*}
avec \begin{enumerate}
\item $\lambda_a=1-\gamma_{a}^{-1}$ pour $a \in \EE$ et $\lambda_a = 1$ pour $a\in\cV_{ell}$ ;
\item $\mu_a=\sum_{s=0}^{\mu-1}\gamma_a^s$ pour $a \in \cV_{ell,\mu}$.
\end{enumerate}
\end{cor}
On en déduit que
$\Delta_0$ admet une présentation en tant que $\ZZ[\Gamma]$-module de la forme
$$\cW_1 \to \cW_0 \to \Delta_0 \to 0 $$
où $\cW_1$ et $\cW_0$ sont des $\ZZ[\Gamma]$-modules libres :
\begin{equation*}
\begin{split}
\cW_0 &= \bigoplus_{a \in \cA} \ZZ[\Gamma] \\
\cW_1 &= \ZZ[\Gamma] \bigoplus \oplus_{a \in \cV_{ell}} \ZZ[\Gamma]
\end{split}
\end{equation*}
La flèche de $\cW_0$ dans $\Delta_0$ est donnée par
$$ (w_a)_{a \in \cA} \mapsto \sum_{a \in \cA} w_a a $$
La flèche de $\cW_1$ dans $\cW_0$ est donnée par
\begin{equation*}
\begin{split}
(v, (w_a)_{a \in \cV_{ell}}) &\mapsto
v (\lambda_a)_{a \in \cA} + (w_{a} \mu_a)_{a \in \cV_{ell}}
\end{split}
\end{equation*}
Plus précisément, le $\ZZ[\Gamma]$-module
$\Delta_0$ admet une résolution projective de la forme
$$ \dots \to \cW_n \to \dots \to \cW_2 \to \cW_1 \to \cW_0 \to \Delta_0 \to 0 $$
où
$\cW_n = \bigoplus_{a \in \cV_{ell}} \ZZ[\Gamma]$ pour $n \geq 2$.
Pour $m \geq 1$, les flèches sont données par
\begin{align*}
\cW_{2m} & \to \cW_{2m-1} &
\cW_{2m+1} & \to \cW_{2m} \\
(w_a)_{a \in\cV_{ell}} & \mapsto (w_a(\gamma_a-1))_{a \in\cV_{ell}} &
(w_a)_{a \in\cV_{ell}} &\mapsto (w_a\mu_a)_{a \in\cV_{ell}}\ .
\end{align*}

\section{Construction d'un \farey normalisé}
\label{dissection}
% !TEX root = ../fareydissection.tex
 % !TEX encoding = IsoLatin
\def\posa{i}
\def\posaast{i^*}
\def\cutun{c_1}
\def\cutdeux{c_2}
\def\kun{k_1}
\def\kdeux{k_2}
\begin{const} L'algorithme décrit dans ce paragraphe permet de construire
à partir d'un \farey $\cF$ de groupe $\Gamma$ un
\farey normalisé $\cF^{norm}$ de même groupe $\Gamma$.
Cette construction est inspirée d'un paragraphe de Siegel dans \cite{siegel},
p. 115 (\S 2) et met en oeuvre l'opération
de dissection sur les surfaces de Riemann (\cite{bost}).
\end{const}
\subsection{Interprétation géométrique}
\label{geometrie}
Nous avons vu qu'à un \farey $\cF$ de groupe $\Gamma$ est associé naturellement
un domaine fondamental $\cU(\cF)$ pour l'action de $\Gamma$ sur $\cH$.
Donnons la description géométrique des opérations de Siegel.
Soit $a$ un élément de $\cV$ tel que $a\neq a^*$, c'est un arc du bord de $\cU(\cF)$.
On choisit deux sommets $s$ et $t$ non elliptiques de $\cU(\cF)$  tels que la géodésique
$a'$ reliant $s$ et $t$ sépare $a$ et $a^*$ ($s$ ou $t$ peut être une extrémité de $a$ ou
de $a^*$). On découpe le domaine le long de $a'$ en deux polygones $\cP_1$ et $\cP_2$
et on translate le polygone contenant $a$ (disons $\cP_1$) par $\gamma_a^{-1}$:
$\cP_1'=\gamma_a^{-1} \cP_1$ : ainsi $a^*$ (resp. $\overline{a}^*$) est un côté du polygone
$\cP_2$ (resp. $\cP'_1$) et l'on peut recoller les deux polygones $\cP_1'$
et $\cP_2$ le long de $a^*$.
On obtient un nouveau domaine fondamental pour $\Gamma$.

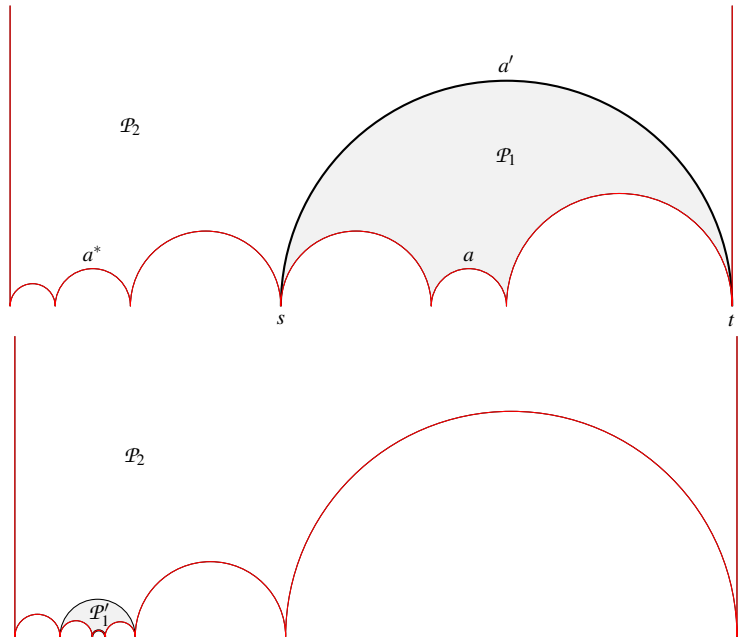
\begin{center}
\begin{tikzpicture} [scale=12]
\scriptsize
\fill [black!5!white] (1/2,0) arc (180:0:1/12) arc (180:0:1/24) arc (180:0:1/8) arc(0:180:1/4);
\draw (1/5,1/3) -- (1/5,0) arc (180:0:1/40)
arc (180:0:1/24) node [midway, above] {$a^*$}
arc (180:0:1/12) arc (180:0:1/12) arc (180:0:1/24)node [midway, above] {$a$} arc (180:0:1/8) -- (1,1/3);
\draw [thick] (1/2,0)  node [below] {$s$} arc(180:0:1/4) node [midway, above] {$a'$} node [below] {$t$};
\draw (1/3,1/5) node {$\cP_2$};
\draw (3/4,1/6) node {$\cP_1$};
\draw[red]  (1/5,1/3) -- (1/5,0) arc (180:0:1/40) arc (180:0:1/24)
            arc (180:0:1/12) arc (180:0:1/12) arc (180:0:1/24) arc (180:0:1/8) -- (1,1/3);
\end{tikzpicture}
\begin{tikzpicture} [scale=12]
\scriptsize
\fill [black!5!white] (1/4,0) arc(180:0:1/56) arc(180:0:1/140) arc(180:0:1/60) arc (0:180:1/24);
\draw (1/5,1/3) -- (1/5,0) arc (180:0:1/40);
\draw (1/4,0) arc (180:0:1/24);
\draw (1/4,0) arc(180:0:1/56);
\draw [thick] (2/7,0) arc(180:0:1/140);
\draw (3/10,0) arc(180:0:1/60) arc (180:0:1/12) arc(180:0:1/4) -- (1,1/3);
\draw (1/3,1/5) node {$\cP_2$};
\draw (5/17,1/40) node {$\cP'_1$};
\draw[red]  (1/5,1/3) -- (1/5,0) arc (180:0:1/40) arc(180:0:1/56) arc(180:0:1/140) arc(180:0:1/60)
                   arc (180:0:1/12) arc(180:0:1/4) -- (1,1/3);
\end{tikzpicture}
\end{center}

Lorsque $a=a^*$, on choisit une coupure $c$ dont l'une des extrémité
est le point elliptique associé à $a$ et l'autre extrémité un point
non elliptique $s$ et la suite du procédé est le même.
Le dessin suivant représente un cas de chemin elliptique d'ordre 3.
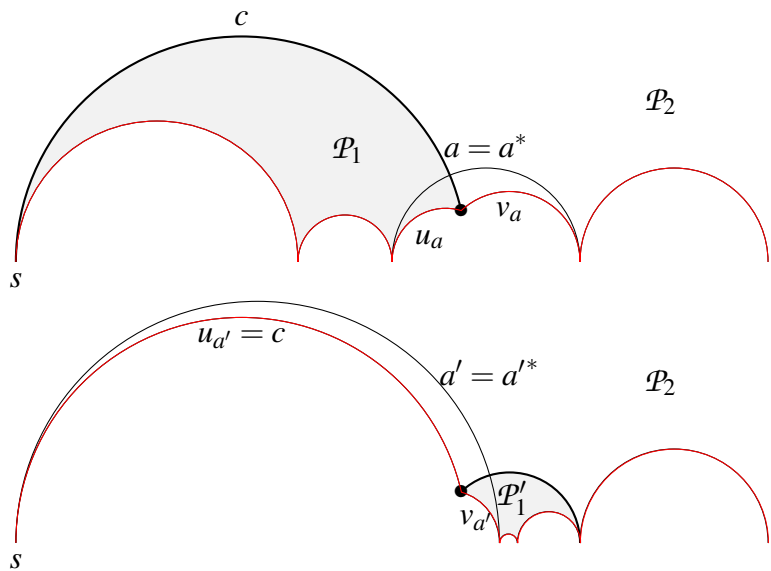
\begin{center}
\begin{tikzpicture}[scale=15]
\fill [black!5!white] (0,0) arc (180:0:1/8) arc (180:0:1/24) arc (180:73.1736:0.0476) arc (13:180:0.2);

\draw (0,0) node [below] {$s$} arc (180:0:1/8) arc (180:0:1/24) arc(180:0:1/12) node [midway, above] {$a=a^*$} arc(180:0:1/12);
\draw (1/3,0) arc (180:73.1736:0.0476) node [midway, below] {$\quad u_a$}node (t) {$\bullet$};
\draw[thick] (0, 0) arc (180:90:0.2) node [above] {$c$} arc (90:13:0.2);
\draw (1/2,0) arc (0:90:0.0625) node [below] {$v_a$} arc (90:133.1736:0.0625);
\draw (4/7,1/7) node {$\cP_2$};
\draw (5/17,1/10) node {$\cP_1$};
\draw[red] (0,0) arc (180:0:1/8) arc (180:0:1/24) arc (180:73.1736:0.0476) arc (133.1736:0:0.0625) arc(180:0:1/12);
\end{tikzpicture}
\begin{tikzpicture}[scale=15]
\fill [black!5!white] (3/7,0) arc (0:74:0.047) arc (133.1736:0:0.0625) arc (0:180:1/36) arc (0:180:1/126) ;
\draw (0,0) node [below] {$s$} arc (180:45:3/14) node [right] {$a'={a'}^*$} arc (45:0:3/14) arc (180:0:1/126) arc (180:0:1/36) arc(180:0:1/12);
\draw (0, 0) arc (180:90:0.2) node [below] {$u_{a'}=c$} arc (90:13:0.2);
\draw (3/7,0) arc (0:74:0.047) node [near end, below] {$v_{a'}$} node {$\bullet$};;
\draw  [thick] (1/2,0) arc (0:133.1736:0.0625);
\draw (4/7,1/7) node {$\cP_2$};
\draw (0.44,0.04) node {$\cP'_1$};
\draw[red] (0,0) arc (180:13:0.2) arc (74:0:0.047) arc (180:0:1/126) arc (180:0:1/36) arc(180:0:1/12);
\end{tikzpicture}
\end{center}
Les deux paragraphes qui suivent décrivent ces opérations en termes de
symboles de Farey uniquement.

Dans la suite, $(\cV, *, \muell)$ désigne un \farey.
On pose $\cV=(a_1, \cdots, a_n)$ et on note  $(\gamma_a)_{a\in \cV}$ 
les données de recollement associées.
Les lettres majuscules désigneront des facteurs de $\cV$, donc formés d'arcs consécutifs.
En particulier, $V$ est le facteur $a_1a_2\cdots a_n$.
Les lettres minuscules réfèrent en général aux arcs de $\cV$.

\subsection{Opération de base (cas non elliptique)}
Choisissons un élément $a$ de $\cV$ tel que $a\neq a^*$ et
un découpage de $V$ de la forme suivante $V=(X_1 a X_2 X_3 a^* X_4)$. On construit
un nouveau \farey $\cF'=(\cV', *, \muell')$
par une opération de base de découpage et recollement de la manière suivante :
on introduit l'arc $a'$ reliant l'extrémité droite du facteur $X_4$
(qui est aussi l'extrémité gauche du facteur $X_1$) à l'extrémité
droite du facteur $X_2$ (qui est aussi l'extrémité gauche du facteur $X_3$)
et on recolle les arcs opposés $a^*$ et $\gamma_a^{-1}(a)$ :

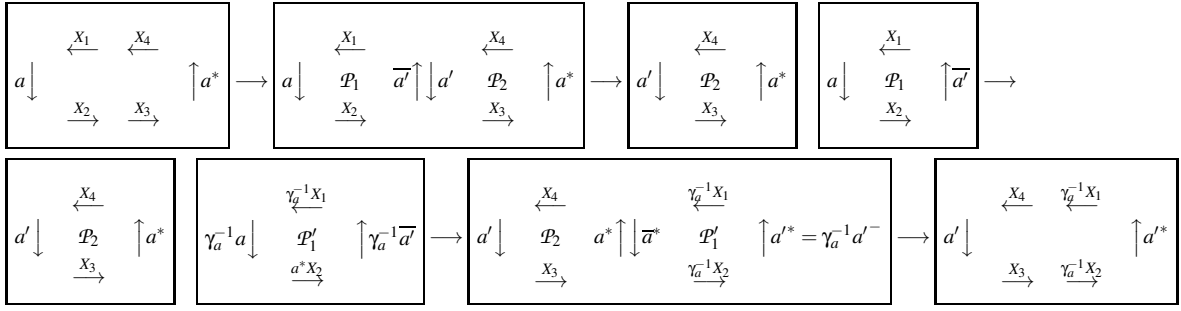
\newcommand\hauteur{{\vline height 1cm depth 0.7cm width 0pt}}
\begin{scriptsize}
\begin{equation*}
\begin{split}
&\boxed{\hauteur\begin{matrix}
&\overset{X_1}{\longleftarrow}&\overset{X_4}{\longleftarrow}\\
a\Big\downarrow&&&\Big\uparrow {a^*}\\
&\overset{X_2}{\longrightarrow}&\overset{X_3}{\longrightarrow}
\end{matrix}}
\longrightarrow
\boxed{\hauteur\begin{matrix}
&\overset{X_1}{\longleftarrow}&&\overset{X_4}{\longleftarrow}\\
a\Big\downarrow&\cP_1&\overline{a'}\Big\uparrow\Big\downarrow {a'} &\cP_2&\Big\uparrow {a^*}\\
&\overset{X_2}{\longrightarrow}&&\overset{X_3}{\longrightarrow}
\end{matrix}}
\longrightarrow
\boxed{\hauteur\begin{matrix}
&\overset{X_4}{\longleftarrow}\\
{a'}\Big\downarrow  &\cP_2&\Big\uparrow {a^*}\\
&\overset{X_3}{\longrightarrow}
\end{matrix}}
\quad
\boxed{\hauteur\begin{matrix}
&\overset{X_1}{\longleftarrow}\\
a\Big\downarrow&\cP_1&\Big\uparrow\overline{a'}\\
&\overset{X_2}{\longrightarrow}
\end{matrix}}
\longrightarrow\\
&
\boxed{\hauteur\begin{matrix}
&\overset{X_4}{\longleftarrow}\\
{a'}\Big\downarrow  &\cP_2&\Big\uparrow {a^*}\\
&\overset{X_3}{\longrightarrow}
\end{matrix}}
\quad
\boxed{\hauteur\begin{matrix}
&\overset{\gamma_a^{-1}X_1}{\longleftarrow}\\
\gamma_a^{-1}a\Big\downarrow &\cP'_1&\Big\uparrow\gamma_a^{-1}\overline{a'}\\
&\overset{a^*X_2}{\longrightarrow}
\end{matrix}}
\longrightarrow
\boxed{\hauteur\begin{matrix}
&\overset{X_4}{\longleftarrow}&&\overset{\gamma_a^{-1} X_1}{\longleftarrow}\\
{a'}\Big\downarrow&\cP_2&{a^*}\Big\uparrow\Big\downarrow \overline{a}^* &\cP_1'&\Big\uparrow {a'}^*={\gamma_a^{-1} {a'}}^{-}\\
&\overset{X_3}{\longrightarrow}&&\overset{\gamma_a^{-1} X_2}{\longrightarrow}
\end{matrix}}
\longrightarrow
\boxed{\hauteur\begin{matrix}
&\overset{X_4}{\longleftarrow}&\overset{\gamma_a^{-1} X_1}{\longleftarrow}\\
{a'}\Big\downarrow&&&\Big\uparrow {a'}^*\\
&\overset{X_3}{\longrightarrow}&\overset{\gamma_a^{-1} X_2}{\longrightarrow}
\end{matrix}}
\end{split}
\end{equation*}
\end{scriptsize}

Cette opération est la traduction en termes de \farey de la description géométrique
sur les domaines fondamentaux du paragraphe \ref{geometrie}.
Décrivons le nouveau \farey $\cF'=(\cV', *, \muell')$. On a
$$V'= (a' X_3\gamma_a^{-1}(X_2) {a'}^* \gamma_a^{-1}(X_1) X_4)$$
Si $v$ est un élément de $\cV$, notons $v'$ l'élément correspondant dans $\cV'$.
Pour $v$ différent de $a$ et $a^*$,
selon la position de $v$,
on a $v'=v$ ou $v'=\gamma_a^{-1} v$.
L'involution $*$ sur $\cF'$ respecte la transformation $v \mapsto v'$,
autrement dit ${v'}^*= {v^*}'$.
Les données de recollement se déduisent de celles de $\cF$.
On a par construction $\gamma_{a'}=\gamma_a$.
Pour $v$ différent de $a$ et $a^*$,
selon les positions respectives de $v$ et de $v^*$,
on a $\gamma_{v'}=\gamma_v$, $\gamma_a^{-1} \gamma_v$, $\gamma_v \gamma_a$ ou
$\gamma_a^{-1} \gamma_v\gamma_a$.
On en déduit que $\gamma_{v'}$ appartient encore à
$\Sl_2(\ZZ)$ (et bien sûr à $\Gamma$)
et donc que le nouveau \farey vérifie la propriété (F4).
Si $v$ est différent de $a$ et de $a^*$, la largeur de $v'$
(déterminant de $A_{v'}$) est égale à la largeur de $v$.
Cela n'est pas vrai pour $a'$ et~${a'}^*={\gamma_a^{-1} {a'}}^{-}$.
Le groupe associé au nouvel \farey $\cF'$ est encore $\Gamma$.

Une transformation similaire est obtenue en échangeant les rôles
de $a$ et de $a^*$ : on choisit alors d'appliquer
$\gamma_a$ au deuxième polygone obtenu par découpage selon
%%%$a'=\gamma_a^{-1}(a'')$,
$a''=\gamma_a^{-1}(a')$:
\renewcommand\hauteur{{\vline height 1.4cm depth 0.7cm width 0pt}}
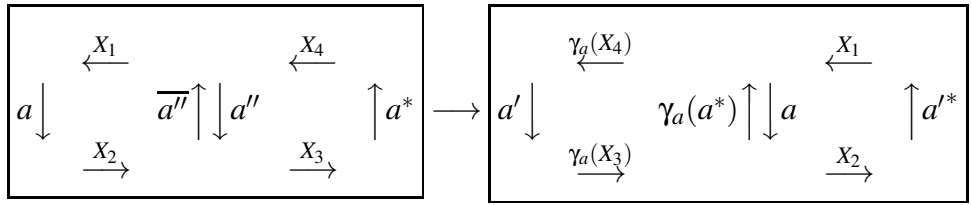
\begin{equation*}
\boxed{\hauteur\begin{matrix}
&\overset{X_1}{\longleftarrow}&&\overset{X_4}{\longleftarrow}\\
a\Big\downarrow&&\overline{a''}\Big\uparrow\Big\downarrow {a''} &&\Big\uparrow {a^*}\\
&\overset{X_2}{\longrightarrow}&&\overset{X_3}{\longrightarrow}
\end{matrix}}
\longrightarrow
\boxed{\hauteur\begin{matrix}
&\overset{\gamma_a(X_4)}{\longleftarrow}&&\overset{X_1}{\longleftarrow}\\
a'\Big\downarrow&&\gamma_a(a^*)\Big\uparrow\Big\downarrow a
&&\Big\uparrow {a'}^*\\
&\overset{\gamma_a(X_3)}{\longrightarrow}&&\overset{X_2}{\longrightarrow}
\end{matrix}}
\end{equation*}

Nous représentons ces opérations de la manière suivante
$$
(\widetilde{\encadre{a} X_2 }\mid X_3 {\encadrebis{a}^*} X_4 \mid \widetilde{X_1})
\mapsto
({a'} X_3\gamma_a^{-1}(X_2){{a'}^*} \gamma_a^{-1} (X_1) X_4)
$$
resp.
$$\left({\encadre{a}} X_2 \mid\widetilde{X_3\encadrebis{a}^* X_4}
\mid{X_1}\right)
\mapsto \left({a'} \gamma_a(X_3)X_2 {a'}^* X_1 \gamma_a(X_4)\right).$$
Les pivots $a$ et $a^*$ que l'on remplace sont soulignés. Les coupures sont indiquées par
un trait vertical. La partie $\widetilde{\cdots}$ entre deux coupures indique les chemins
transformés par $\gamma_a^{-1}$ (resp. $\gamma_a$).
De manière équivalente, les opérations peuvent être représentées par
$$(\widetilde{X_1\encadre{a} X_2 }\mid X_3 {\encadrebis{a}^*} X_4 \mid)
\mapsto ( X_4 {a'} X_3 \gamma_a^{-1}(X_2){{a'}^*} \gamma_a^{-1}(X_1))$$
resp.
$$(X_1\encadre{a} X_2 \mid \widetilde{X_3{\encadrebis{a}^*} X_4}\mid)
\mapsto (\gamma_a(X_4) {a'} \gamma_a(X_3) X_2{{a'}^*} X_1).$$

\subsection{Opération de base (cas elliptique)}

Supposons maintenant que $a^*=a$. Les opérations de base sont les suivantes
\begin{equation*}
(\widetilde{X_1 \encadre{a}} \mid X_2\mid)\rightarrow (a' X_1 \gamma_a(X_2))
\quad \text{resp.} \quad
(X_1 \encadre{a}\mid \widetilde{X_2}\mid)\rightarrow (a' \gamma_a^{-1}(X_1)X_2)
\ .
\end{equation*}
Pour comprendre l'opération, il faut se souvenir que le polygone fondamental construit
à partir de $\cF$ ne contient pas $a=(r_a,s_a)$ comme bord, mais les arcs
$u_a=(r_a,t_a)$ et $v_a=(t_a,s_a)$ où $t_a$ appartient au demi-plan de Poincaré $\cH$.
La coupure $c=(s,t_a)$ est faite entre $t_a$ et le sommet $s$ extrémité commune de
$X_1$ et $X_2$.
La dissection faite sur ce polygone est alors (dans le premier cas)
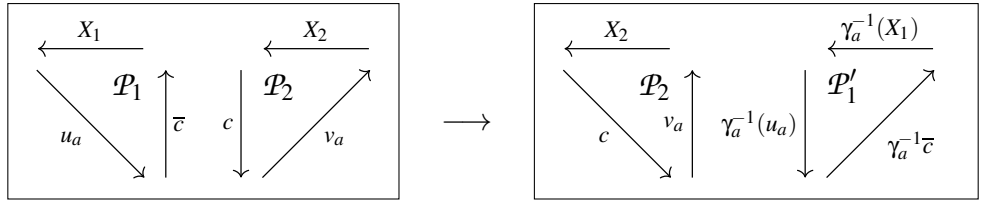
\begin{center}
\begin{tikzpicture}[commutative diagrams/every diagram]
\draw (-0.1,0) rectangle (5.1,2.6);
\draw (6.9,0) rectangle (12.8,2.6);
\draw (6,1) node {$\longrightarrow$};
\node (P1) at (0,2){\phantom{X}};
\node (P2) at (2,2){\phantom{X}};
\node (P3) at (3,2){\phantom{X}};
\node (P4) at (5,2){\phantom{X}};
\node (P5) at (7,2){\phantom{X}};
\node (M1) at (1.5,1.5){$\cP_1$};
\node (M2) at (8.5,1.5){$\cP_2$};
\node (P6) at (9,2){\phantom{X}};
\node (P7) at (10.5,2){\phantom{X}};
\node (P8) at (12.5,2){\phantom{X}};
\node (Q1) at (2,0){\phantom{X}};
\node (Q2) at ((3,0){\phantom{X}};
\node (Q3) at (9,0){\phantom{X}};
\node (Q4) at (10.5,0){\phantom{X}};
\node (M3) at (3.5,1.5){$\cP_2$};
\node (M4) at (11,1.5){$\cP_1'$};

\path[commutative diagrams/.cd, every arrow, every label]
(P2) edge node [swap] {$X_1$} (P1)
(P1) edge node [swap] {$u_a$} (Q1)
(Q1) edge node [swap] {$\overline{c}$} (P2)
(P4) edge node [swap] {$X_2$} (P3)
(P3) edge node [swap] {$c$} (Q2)
(Q2) edge node [swap] {$v_a$} (P4)
(P6) edge node [swap] {$X_2$} (P5)
(P5) edge node [swap] {$c$} (Q3)
(Q3) edge node {$v_a$} (P6)
(P8) edge node [swap] {$\gamma_a^{-1}(X_1)$} (P7)
(P7) edge node [swap] {$\gamma_a^{-1}(u_a)$} (Q4)
(Q4) edge node [swap] {$\gamma_a^{-1}\overline{c}$} (P8);
\end{tikzpicture}
\end{center}
On a $\gamma_a^{-1}(u_a)=\overline{v_a}$ et si $a'=c c^*$ avec $c^*=\gamma_a^{-1}\overline{c}$,
on a $\gamma_{a'}=\gamma_a$. Le second cas se traite
de la même manière en appliquant $\gamma_a$ à $\cP_2$.

\subsection{Étape de Siegel}
Notons $W$ un facteur initial de $V$ qui est normalisé
(c'est-à-dire lui-même formé de facteurs de la forme
$ab{a}^*{b}^*$, $aa^*$ ou $a$ avec $a=a^*$)
et posons $V=(WX)$.
L'étape de Siegel que nous allons décrire
permet de remplacer le \farey initial $\cF$ par
un \farey $\cF'$ avec $V'=(W'X')$ et $W'$ strictement plus long
que $W$ en (au plus) 4 opérations de base.

%%L'étape de Siegel permet d'augmenter la longueur du segment «initial» $W$
%%du symbole de Farey
%%qui est normalisé (c'est-à-dire formé de facteurs de la forme
%%$ab{a}^*{b}^*$, $aa^*$ ou $a$ avec $a=a^*$)
%%en (au plus) 4 opérations de base.
Dans la suite, on omettra les $'$ dans la notation : par exemple,
$$(\widetilde{X_1\encadre{a} X_2 }\mid X_3 {\encadre{a}^*} X_4 \mid)
\mapsto ( X_4 {a'} X_3 \gamma_a(X_2){{a'}^*} \gamma_a(X_1))=
( X_4' {a'} X_3' X_2'{a'}^* X_1')$$ sera simplement écrit
$$(\widetilde{X_1\encadre{a} X_2 }\mid X_3{\encadre{a}^*} X_4 \mid)
\mapsto ( X_4 {a} X_3 X_2a^*X_1)\ .$$
En pratique, nous supposerons que $W$ contient l'arc $(\infty,0)$.
Le choix de la coupure qui sera modifiée par un élément de $\Gamma$
est fait de manière à ne pas transformer l'arc $\{\infty,0\}$,
ce qui revient à ne pas appliquer un élément de $\Gamma$ sur $W$.
Expérimentalement, cela diminue de manière drastique la taille
des coefficients des matrices et le temps de calcul.

%À une rotation près, on peut supposer que la partie normalisée $W$ de
%$V$ est au début de $V$.

\begin{enumerate}
\item Cas parabolique. On considère un arc $a$ à distance 1 de
$a^*$ (appelé pivot). On peut supposer qu'il précède $a^*$ quitte à changer de notation.
%soit l'arc suivant $a$, autrement dit, $\kun^*=\kun+1$.
La longueur de $W$ augmente de 2 avec l'une des deux opérations suivantes.

\begin{equation*}
\begin{split}
(WXa{a^*}Y)=(WX\encadre{a}\mid {\encadrebis{a^*}}\widetilde{Y}\mid ) \rightarrow (a{a^*}WXY)
%\text{ i.e } [\posa,\cutun,\cutdeux]=[\kun, \kun^*,w_1]=[\kun, \kun+1,1]
\end{split}
\end{equation*}
ou
\begin{equation*}
\begin{split}
(WXa{a^*}Y)=(W \mid \widetilde{X}\encadre{a}\mid {\encadrebis{a^*}}Y) \rightarrow (a{a^*}XY W)
\overset{rotation}{\longrightarrow} (W a{a^*} XY)\ .
%\text{ i.e } [\posa,\cutun,\cutdeux]=[\kun, \kun^*=\kun+1,w_1]=[\kun, \kun+1,lg(W)]
\end{split}
\end{equation*}

\item Cas hyperbolique.
On considère deux arcs $a$ et $b$ liés avec $a$ précédant $b$ dans
$V$ (appelés pivots de l'étape).
%de position respective $\kun$ et $\kdeux$. On a donc $\kun < \kdeux < \kun^* < \kdeux$.
Plusieurs suites d'opérations sont possibles.
Quand il y a plusieurs choix possibles, les coupures sont choisies
de manière à ce que le pivot non utilisé soit dans le
même facteur que son image par l'involution.
Dans tous les cas, la longueur de $W$ augmente de 4.
En pratique, nous n'avons utilisé que le premier cas.
\begin{enumerate}
\item $b$ suit $a$.
%%Cas où $\kdeux=\kun+1$.
\begin{equation*}
\begin{split}
&(WXabY{a^*}Z{b^*}T)=(WXa\encadre{b}Y\mid\widetilde{{a^*}Z \encadrebis{b^*}T}\mid)
%%\text{ i.e } [\posa,\cutun,\cutdeux]=[\kdeux,\kun^*,w_1]
\\
&\rightarrow (b{a^*}ZY{b^*}WXaT)=(\widetilde{b\encadrebis{a^*}ZY}\mid {b^*}WX\encadre{a}T\mid)
%%\text{ i.e } [\posa,\cutun,\cutdeux]=[\kun,\kdeux,\kdeux^*]
\\
&\rightarrow (a^*b^* W XZYa b T)=(a^*\encadre{b^*} W XZYa \mid \widetilde{\encadrebis{b} T}\mid)
%%\text{ i.e } [\posa,\cutun,\cutdeux]=[\kdeux^*,\kdeux,\kun^*]
\\
&\rightarrow ({b^*}WXZYa b a^*T)=({b^*}WXZY \encadre{a}\mid\widetilde{ b \encadrebis{a^*}T}\mid)
%%\text{ i.e } [\posa,\cutun,\cutdeux]=[\kun,\kun^*,\kdeux^*]
%%ou [\posa,\cutun,\cutdeux]=[\kun,k_2,\kdeux^*]
\\
&\rightarrow (ab{a^*}{b^*} WXZYT)
\end{split}
\end{equation*}

%%\begin{center}
%%\begin{tabular}{|c|c|c|c|c|}
%%\hline &X &Y &Z &T\\\hline
%%X& 0&1 &2 &3 \\\hline
%%Y& 1&2 &3 &4\\\hline
%%Z&2 &3 &4 &5\\\hline
%%T&3 &4 &5 &6\\\hline
%%\end{tabular}
%%\end{center}
\item $b$ suit $a$ (autre méthode).
\begin{equation*}
\begin{split}
&(WXabY{a^*}Z{b^*}T)=(W\mid \widetilde{Xa\encadre{b}Y}\mid{a^*}Z \encadrebis{b^*}T)
%%\text{ i.e } [\posa,\cutun,\cutdeux]=[\kdeux,\kun^*,w_1]
\\
&\rightarrow (b{a^*}ZY{b^*}XaTW)=(\widetilde{b\encadrebis{a^*}ZY}\mid {b^*}X\encadre{a}TW\mid)
%%\text{ i.e } [\posa,\cutun,\cutdeux]=[\kun,\kdeux,\kdeux^*]
\\
&\rightarrow (a^*b^* XZYa b TW)=(\widetilde{a^*\encadre{b^*}}\mid XZYa \encadrebis{b} T W \mid)
%%\text{ i.e } [\posa,\cutun,\cutdeux]=[\kdeux^*,\kdeux,\kun^*]
\\
&\rightarrow ({b^*}XZYa b a^*TW)=(\widetilde{{b^*}XZY \encadre{a}b} \mid \encadrebis{a^*}TW\mid)
%%\text{ i.e } [\posa,\cutun,\cutdeux]=[\kun,\kun^*,\kdeux^*]
%ou [\posa,\cutun,\cutdeux]=[\kun,k_2,\kdeux^*]
\\
&\rightarrow (ab{a^*}{b^*}XZYTW) \overset{rotation}{\longrightarrow}(Wab{a^*}{b^*}XZYT)
\end{split}
\end{equation*}
%%%\begin{center}
%%%\begin{tabular}{|c|c|c|c|c|}
%%%\hline &X &Y &Z &T\\\hline
%%%X& 4&5 &4 &2 \\\hline
%%%Y& 5 &6 &5&3\\\hline
%%%Z& 4 &5 &4 &2\\\hline
%%%T& 2 &3 &2 &0\\\hline
%%%\end{tabular}
%%%\end{center}

\item Cas général.
\begin{equation*}
\begin{split}
&(WXaUbY{a^*}Z{b^*}T)=(WX\encadre{a}Ub\mid \widetilde{Y{\encadrebis{a^*}}Z}\mid {b^*}T)
%%\text{ i.e } [\posa,\cutun,\cutdeux]=[\kun,\kdeux+1, \kdeux^*]
\\
&\rightarrow (aYUb{a^*}{b^*}TWXZ)=(aYU\encadre{b}\mid \widetilde{{a^*}{\encadrebis{b^*}}T}\mid WXZ)
%%\text{ i.e } [\posa,\cutun,\cutdeux]=[\kdeux,\kun^*, w_1]
\\
&\rightarrow (b{a^*}{b^*}WXZaYUT)=(b\encadre{{a^*}}{b^*}W\mid \widetilde{XZ\encadrebis{a}YUT}\mid )
%%\text{ i.e } [\posa,\cutun,\cutdeux]=[\kun^*,w_1+n, \kdeux]
\\
&\rightarrow ({a^*}XZ {b^*}WabYUT)=(\widetilde{{a^*}XZ \encadre{b^*}}\mid Wa\encadrebis{b}YUT\mid )
%%\text{ i.e } [\posa,\cutun,\cutdeux]=[\kdeux^*,w_1, \kun^*]
\\
&\rightarrow ({b^*}Wab{a^*}XZYUT)=(\widetilde{{b^*}W\encadre{a}b}\mid {\encadrebis{a^*}}XZYUT\mid )
%%\text{ i.e } [\posa,\cutun,\cutdeux]=[\kun, \kun^*, \kdeux^*]
\\
&\rightarrow (ab{a^*}{b^*} WXZYUT)
\end{split}
\end{equation*}
%%%\begin{center}\begin{tabular}{|c|c|c|c|c|c|}
%%%\hline &X &U &Y &Z &T\\\hline
%%%X&4&3&4&3&5\\\hline
%%%U&3&2&3&4&3\\\hline
%%%Y&4&4&4&5&4\\\hline
%%%Z&3&5&6&5&4\\\hline
%%%T&5&4&5&4&3\\\hline
%%%\end{tabular}
%%%\end{center}
\item Cas général avec moins d'étapes.
\begin{equation*}
\begin{split}
&(WXaUbY{a^*}Z{b^*}T)=(W\mid \widetilde{X\encadre{a}}\mid UbY\encadrebis{a^*}Z b^*T)
%\text{ i.e } [\posa,\cutun,\cutdeux]=[\kun,\kun+1,w_2+1]
\\
&\rightarrow (W a U b Y a^* X Z b^*T)=(W\mid \widetilde{a U \encadre{b} Y a^* }\mid X Z \encadrebis{b^*} T)
%\text{ i.e } [\posa,\cutun,\cutdeux]=[\kdeux,\kun^*+1,\kun]
\\
&\rightarrow (T W b X Z Y a^* b^* a U)=(T W b \mid \widetilde{X Z Y \encadre{a^*}} \mid b^* \encadrebis{a} U)
%\text{ i.e } [\posa,\cutun,\cutdeux]=[\kun^*,\kdeux^*,\kdeux+1]
\\
&\rightarrow (UTW b a b^*a^* XZY)
\overset{rotation}{\longrightarrow} (Wb a b^*a^* XZYUW)
\end{split}
\end{equation*}

%%%\begin{center}\begin{tabular}{|c|c|c|c|c|c|}
%%%\hline&X &U &Y &Z &T \\\hline
%%%X&4 &3&4 &3 &2 \\\hline
%%%U&3 &2 &3 &2 &1\\\hline
%%%Y&4 &3 &4 &3 &2 \\\hline
%%%Z&3 &2 &3 &2 &1 \\\hline
%%%T&2 &1 &2 &1 &0\\\hline
%%%
%%%\end{tabular}
%%%\end{center}
%%\begin{equation*}
%%\begin{split}
%%&(WaUb{a^*}Z{b^*}T)=(W\mid \widetilde{a U \encadre{b} a^* }\mid Z \encadrebis{b^*} T)
%%%%\text{ i.e } [\posa,\cutun,\cutdeux]=[\kun,\kdeux,\kdeux^*]
%%\\
%%&\rightarrow (T W b Z a^* b^* a U)=(T W b \mid \widetilde{ Z \encadre{a^*}} \mid b^* \encadre{a} U)
%%%%\text{ i.e } [\posa,\cutun,\cutdeux]=[\kdeux^*,\kdeux,\kun^*]
%%\\
%%&\rightarrow (UTW b a b^*a^* Z)
%%\end{split}
%%\end{equation*}

\end{enumerate}

\item Cas elliptique. Le pivot est ici un élément $a$ de $\cV_{ell}$.
%On a $\kun^*=\kun$ et on fait l'opération
\begin{equation*}
\begin{split}
(\mid WX\encadre{a}\widetilde{Y})\rightarrow (WXYa)
 \overset{rotation}{\longrightarrow} (aWXY)
%\text{ i.e } [\posa,\cutun,\cutdeux]=[\kun, k_1^*,w_1]
\end{split}
\end{equation*}
Deuxième possibilité :
\begin{equation*}
\begin{split}
( W\mid \widetilde{X}\encadre{a}Y)\rightarrow (aWXY)
\end{split}
\end{equation*}

La longueur de $W$ augmente de 1.
\end{enumerate}

\subsection{Algorithme complet}
En alternant ces différentes opérations, on obtient l'algorithme suivant.
L'algorithme modifie en place un \farey $(\cV,*,\muell)$ par une suite
d'opérations de Siegel. Comme dans le paragraphe précédent, on note $V$ le
facteur $a_1a_2\dots a_n$ de $\cV = (a_1,\dots,a_n)$ et $W$ est un préfixe de
$V$ normalisé, formé de facteurs de la forme $aba^*b^*$, $aa^*$ ou $a$
avec $a = a^*$.
\begin{algorithm}[H]
\caption{Normalisation d'un \farey}
\begin{algorithmic}[1]
\REQUIRE Un \farey $(\cV,*,\muell)$.
\ENSURE Un \farey normalisé
\STATE $W \gets \emptyset$ \COMMENT{le facteur vide}
\WHILE{$|W| < |V|$}
  \IF{il existe un facteur $W'$ de $V$ normalisé contenant strictement $W$ }
    \STATE $W \gets W'$;
    \COMMENT {$|W|$ augmente de $|W'| - |W|$.}
  \ELSIF{il existe $a$ hors de $W$ tel que $d(a,a^*)=0$}
    \STATE appliquer l'opération de Siegel elliptique de pivot $a$;
    \COMMENT {$|W|$ augmente de $1$.}
  \ELSIF{il existe $a$ hors de $W$ tel que $d(a,a^*)=1$}
    \STATE appliquer l'opération de Siegel parabolique de pivot $a$;
    \COMMENT {$|W|$ augmente de $2$.}
 \ELSIF{il existe hors de $W$ deux arcs entrelacés $a$ et $b$}
    \STATE appliquer l'opération de Siegel hyperbolique de pivots $a$ et $b$;
    \COMMENT {$|W|$ augmente de $4$.}
   \ENDIF
\ENDWHILE
\RETURN le \farey obtenu.
\end{algorithmic}
\end{algorithm}

\begin{proof}
Pour voir que l'algorithme termine, il suffit de remarquer que tant que $W$
est différent de $V$, au moins l'une des quatre conditions a lieu. En effet,
si aucune des trois premières conditions n'est vérifiée, on cherche le
premier arc $f$ suivant $W$ tel que $f^*$ précède $f$. Ni $f$, ni $f^*$ ne
sont dans $W$. On prend alors comme pivot $a=f^*$ et $b$ l'arc suivant $a$,
les arcs $a$ et $b$ étant alors nécessairement entrelacés.
\end{proof}

\subsection{Un exemple détaillé : $\Gamma_0(22)$}

L'entrée est un symbole de Farey arbitraire pour $\Gamma_0(22)$, ici un symbole
de Farey unimodulaire obtenu avec l'algorithme de Kulkarni. La partie $W$ déjà
normalisée est indiquée entre parenthèses, les deux pivots choisis pour
l'étape sont soulignés:
\smallskip

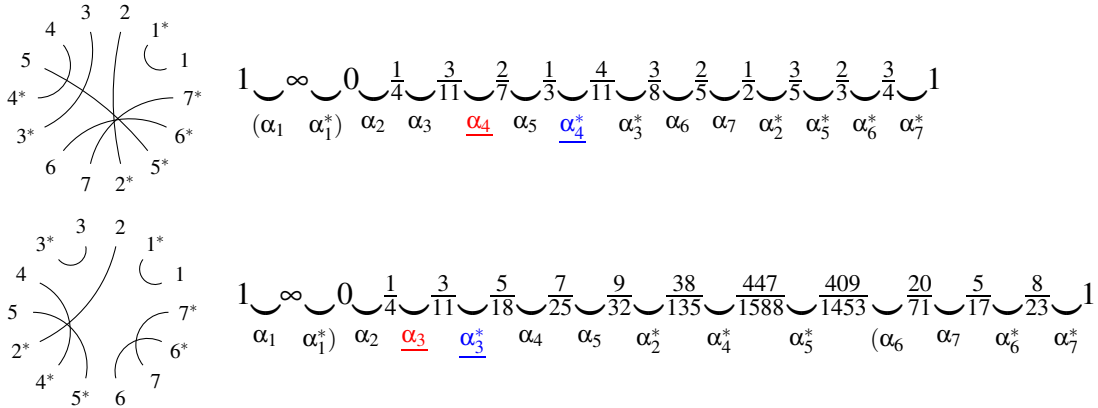
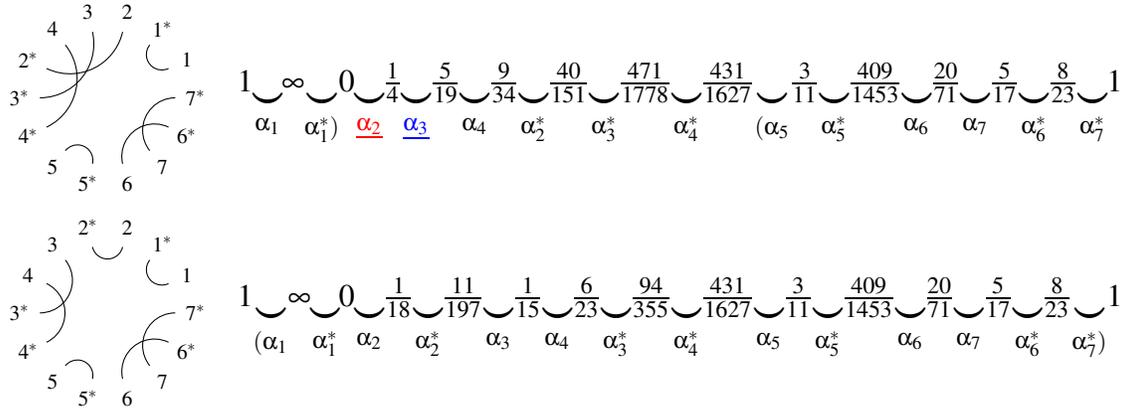
\begin{tabular}{m{2.7cm}m{10cm}}
\begin{tikzpicture}[scale=0.9]
\scriptsize
\draw \link {1}{$1$}{2}{$1^*$}{25.7143};
\draw \link {3}{$2$}{11}{$2^*$}{25.7143};
\draw \link {4}{$3$}{8}{$3^*$}{25.7143};
\draw \link {5}{$4$}{7}{$4^*$}{25.7143};
\draw \link {6}{$5$}{12}{$5^*$}{25.7143};
\draw \link {9}{$6$}{13}{$6^*$}{25.7143};
\draw \link {10}{$7$}{14}{$7^*$}{25.7143};
\end{tikzpicture}
&
$1
\underset{ {(\aaa_{1}}}{\underparen{}}\infty
\underset{ {\aaa^*_{1})}}{\underparen{}}0
\underset{ {\aaa_{2}}}{\underparen{}}\frac{1}{4}
\underset{ {\aaa_{3}}}{\underparen{}}\frac{3}{11}
\underset{ {\encadre{\aaa_{4}}}}{\underparen{}}\frac{2}{7}
\underset{ {\aaa_{5}}}{\underparen{}}\frac{1}{3}
\underset{ {\encadrebis{\aaa^*_{4}}}}{\underparen{}}\frac{4}{11}
\underset{ {\aaa^*_{3}}}{\underparen{}}\frac{3}{8}
\underset{ {\aaa_{6}}}{\underparen{}}\frac{2}{5}
\underset{ {\aaa_{7}}}{\underparen{}}\frac{1}{2}
\underset{ {\aaa^*_{2}}}{\underparen{}}\frac{3}{5}
\underset{ {\aaa^*_{5}}}{\underparen{}}\frac{2}{3}
\underset{ {\aaa^*_{6}}}{\underparen{}}\frac{3}{4}
\underset{ {\aaa^*_{7}}}{\underparen{}}1
 $ \\
\end{tabular}

\begin{tabular}{m{2.7cm}m{10cm}}
\begin{tikzpicture}[scale=0.9]
\scriptsize
\draw \link {1}{$1$}{2}{$1^*$}{25.7143};
\draw \link {3}{$2$}{8}{$2^*$}{25.7143};
\draw \link {4}{$3$}{5}{$3^*$}{25.7143};
\draw \link {6}{$4$}{9}{$4^*$}{25.7143};
\draw \link {7}{$5$}{10}{$5^*$}{25.7143};
\draw \link {11}{$6$}{13}{$6^*$}{25.7143};
\draw \link {12}{$7$}{14}{$7^*$}{25.7143};
\end{tikzpicture}
&
$ 1
\underset{ {\aaa_{1}}}{\underparen{}}\infty
\underset{ {\aaa^*_{1})}}{\underparen{}}0
\underset{ {\aaa_{2}}}{\underparen{}}\frac{1}{4}
\underset{ {\encadre{\aaa_{3}}}}{\underparen{}}\frac{3}{11}
\underset{ {\encadrebis{\aaa^*_{3}}}}{\underparen{}}\frac{5}{18}
\underset{ {\aaa_{4}}}{\underparen{}}\frac{7}{25}
\underset{ {\aaa_{5}}}{\underparen{}}\frac{9}{32}
\underset{ {\aaa^*_{2}}}{\underparen{}}\frac{38}{135}
\underset{ {\aaa^*_{4}}}{\underparen{}}\frac{447}{1588}
\underset{ {\aaa^*_{5}}}{\underparen{}}\frac{409}{1453}
\underset{ {(\aaa_{6}}}{\underparen{}}\frac{20}{71}
\underset{ {\aaa_{7}}}{\underparen{}}\frac{5}{17}
\underset{ {\aaa^*_{6}}}{\underparen{}}\frac{8}{23}
\underset{ {\aaa^*_{7}}}{\underparen{}}1
 $\\
\end{tabular}

\begin{tabular}{m{2.7cm}m{10cm}}
\begin{tikzpicture}[scale=0.9]
\scriptsize
\draw \link {1}{$1$}{2}{$1^*$}{25.7143};
\draw \link {3}{$2$}{6}{$2^*$}{25.7143};
\draw \link {4}{$3$}{7}{$3^*$}{25.7143};
\draw \link {5}{$4$}{8}{$4^*$}{25.7143};
\draw \link {9}{$5$}{10}{$5^*$}{25.7143};
\draw \link {11}{$6$}{13}{$6^*$}{25.7143};
\draw \link {12}{$7$}{14}{$7^*$}{25.7143};
\end{tikzpicture}
&
$ 1
\underset{ {\aaa_{1}}}{\underparen{}}\infty
\underset{ {\aaa^*_{1})}}{\underparen{}}0
\underset{ {\encadre{\aaa_{2}}}}{\underparen{}}\frac{1}{4}
\underset{ {\encadrebis{\aaa_{3}}}}{\underparen{}}\frac{5}{19}
\underset{ {\aaa_{4}}}{\underparen{}}\frac{9}{34}
\underset{ {\aaa^*_{2}}}{\underparen{}}\frac{40}{151}
\underset{ {\aaa^*_{3}}}{\underparen{}}\frac{471}{1778}
\underset{ {\aaa^*_{4}}}{\underparen{}}\frac{431}{1627}
\underset{ {(\aaa_{5}}}{\underparen{}}\frac{3}{11}
\underset{ {\aaa^*_{5}}}{\underparen{}}\frac{409}{1453}
\underset{ {\aaa_{6}}}{\underparen{}}\frac{20}{71}
\underset{ {\aaa_{7}}}{\underparen{}}\frac{5}{17}
\underset{ {\aaa^*_{6}}}{\underparen{}}\frac{8}{23}
\underset{ {\aaa^*_{7}}}{\underparen{}}1
 $\\
\end{tabular}

\begin{tabular}{m{2.7cm}m{10cm}}
\begin{tikzpicture}[scale=0.9]
\scriptsize
\draw \link {1}{$1$}{2}{$1^*$}{25.7143};
\draw \link {3}{$2$}{4}{$2^*$}{25.7143};
\draw \link {5}{$3$}{7}{$3^*$}{25.7143};
\draw \link {6}{$4$}{8}{$4^*$}{25.7143};
\draw \link {9}{$5$}{10}{$5^*$}{25.7143};
\draw \link {11}{$6$}{13}{$6^*$}{25.7143};
\draw \link {12}{$7$}{14}{$7^*$}{25.7143};
\end{tikzpicture}
&
$ 1
\underset{ {(\aaa_{1}}}{\underparen{}}\infty
\underset{ {\aaa^*_{1}}}{\underparen{}}0
\underset{ {\aaa_{2}}}{\underparen{}}\frac{1}{18}
\underset{ {\aaa^*_{2}}}{\underparen{}}\frac{11}{197}
\underset{ {\aaa_{3}}}{\underparen{}}\frac{1}{15}
\underset{ {\aaa_{4}}}{\underparen{}}\frac{6}{23}
\underset{ {\aaa^*_{3}}}{\underparen{}}\frac{94}{355}
\underset{ {\aaa^*_{4}}}{\underparen{}}\frac{431}{1627}
\underset{ {\aaa_{5}}}{\underparen{}}\frac{3}{11}
\underset{ {\aaa^*_{5}}}{\underparen{}}\frac{409}{1453}
\underset{ {\aaa_{6}}}{\underparen{}}\frac{20}{71}
\underset{ {\aaa_{7}}}{\underparen{}}\frac{5}{17}
\underset{ {\aaa^*_{6}}}{\underparen{}}\frac{8}{23}
\underset{ {\aaa^*_{7})}}{\underparen{}}1
 $\\
\end{tabular}

Ce dernier symbole est normalisé.

  %STOP

\subsection{Conclusion}

La normalisation d'un symbole de Farey se révèle très coûteuse: la taille des
rationnels $p/q$ intervenant dans le symbole normalisé augmente
exponentiellement dans notre implantation pour $\Gamma = \Gamma_0(N)$. Plus
précisément, la hauteur naïve $\log_2 \max(|p|,|q|)$ des rationnels
intervenant dans le symbole fourni par l'algorithme de Kulkarni est $O(\log
N)$. Nous ne savons pas borner efficacement la hauteur des rationnels
intervenant dans sa normalisation : elle est trivialement en $O(2^N)$;
expérimentalement, elle se révèle un peu inférieure à $N/2$ pour $N \leq
40000$. Notre implantation de l'algorithme de Kulkarni calcule un symbole de
Farey unimodulaire pour $\Gamma_0(40000)$ en 6 secondes, et il faut environ
une heure pour le normaliser.

Il serait intéressant d'estimer finement la taille du symbole de Farey
normalisé obtenu, mais surtout de la minimiser, en changeant l'ordre des
opérations ou en utilisant d'autres opérations de base. Nous n'y sommes pas
parvenus. Nous n'avons pas non plus d'application algorithmique dans laquelle
utiliser un symbole normalisé serait plus avantageux qu'utiliser un symbole
de Farey unimodulaire. Nous sommes donc intéressés par toute suggestion.

\appendix
\label{fareyexemple}
% !TEX root = ../fareydissection.tex
% !TEX encoding = IsoLatin

\def\bsmall{\begingroup\scriptsize}
\def\esmall{\endgroup}

\section{Exemples}

Donnons quelques exemples que l'on peut déduire des implémentations dans
Pari/GP. Le groupe $\Gamma$ est toujours le sous-groupe de congruence
$\Gamma_0(N)$. Dans le cas où le symbole de Farey unimodulaire n'est pas normalisé,
nous donnons aussi un symbole de Farey normalisé. La syntaxe GP est
\begin{verbatim}
mspolygon(N,0) \\ symbole de Farey unimodulaire
mspolygon(N,1) \\ symbole de Farey normalisé
\end{verbatim}

Nous donnons le symbole de Farey, la représentation de l'involution
du symbole de Farey unimodulaire et du symbole de Farey normalisé et dans les petits cas les
polygones hyperboliques obtenus.
Les noms $\aaa_i$ permettent de donner l'involution (ils sont bien sûr différents
pour le symbole de Farey unimodulaire et le symbole de Farey normalisé).

%%\begin{equation*}
%%\Gamma_0(1) :
%%\input fareydraw/0_1.tex
%%\end{equation*}

%\begin{center}
%\begin{tikzpicture}
%\input fareydraw/cercle0_1.tex
%\end{tikzpicture}
%\end{center}

%%\begin{center}
%%\input polygondraw/2_1.tex
%%\end{center}

\vskip 1cm

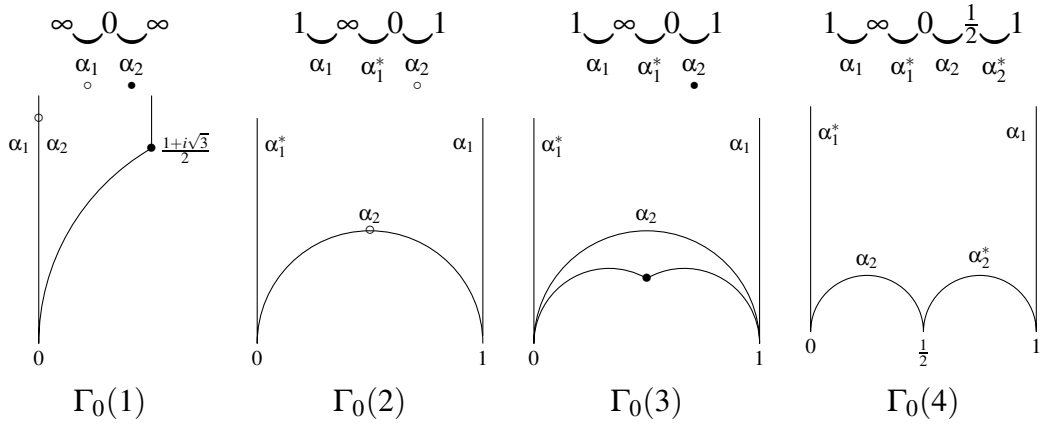
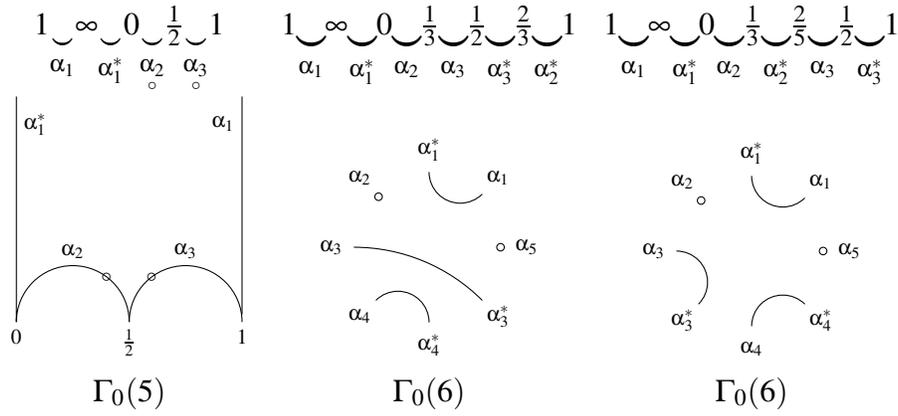
$$
\begin{matrix}
\infty\underset{\underset{\circ}{\aaa_{1}}}{\underparen{}}0
\underset{\underset{\bullet}{\aaa_{2}}}{\underparen{}}\infty
&
1
\underset{ {\aaa_{1}}}{\underparen{}}\infty
\underset{ {\aaa^*_{1}}}{\underparen{}}0
\underset{\underset{\circ}{\aaa_{2}}}{\underparen{}}1

&
1
\underset{ {\aaa_{1}}}{\underparen{}}\infty
\underset{ {\aaa^*_{1}}}{\underparen{}}0
\underset{\underset{\bullet}{\aaa_{2}}}{\underparen{}}1

&
1
\underset{ {\aaa_{1}}}{\underparen{}}\infty
\underset{ {\aaa^*_{1}}}{\underparen{}}0
\underset{ {\aaa_{2}}}{\underparen{}}\frac{1}{2}
\underset{ {\aaa^*_{2}}}{\underparen{}}1

\\
\begin{tikzpicture}[scale=3]
\scriptsize
\draw  (0,1.1) -- (0,1) node {$\circ$} -- (0,0) node [very near start, left] {$\aaa_1$} node [very near start, right] {$\aaa_2$}
node [below] {$0$}arc (180:120:1) node {$\bullet$} node [right] {$\frac{1+i\sqrt3}2$} -- (1/2,1.1);
    \end{tikzpicture}
&
\begin{tikzpicture}[scale=3]
\scriptsize
\draw (0,1)--(0,0) node [very near start, right] {$\aaa_1^*$} node [below] {$0$}node (start) {} arc (180:0:1/2)  node [midway, above] {$\aaa_{2}$} node (start) {} (1/2,1/2) node {$\circ$};
 \draw (start) node{} node [below] {$1$} -- (1,1) node [very near end, left] {$\aaa_1$};
    \end{tikzpicture}
&
\begin{tikzpicture}[scale=3]
\scriptsize
\draw (0,1)--(0,0) node [very near start, right] {$\aaa_1^*$} node [below] {$0$}node (start) {} arc (180:0:1/2)  node [midway, above] {$\aaa_{2}$};
 \draw (start) arc (180:60:0.33333333) node {$\bullet$} arc (120:0:0.33333333) node [below] {$1$} -- (1,1) node [very near end, left] {$\aaa_1$};
 \end{tikzpicture}
&
\begin{tikzpicture}[scale=3]
\scriptsize
\draw (0,1)--(0,0) node [very near start, right] {$\aaa_1^*$} node [below] {$0$}node (start) {} arc (180:0:0.25000000)  node [midway, above] {$\aaa_{2}$} node [below]{$\frac12$}
node (start) {} arc (180:0:1/4)  node [midway, above] {$\aaa^*_{2}$} node [below] {$1$} -- (1,1) node [very near end, left] {$\aaa_1$};
 \end{tikzpicture}
 \\
 \Gamma_0(1)&\Gamma_0(2)&\Gamma_0(3)&\Gamma_0(4)
\end{matrix}
$$

\vskip 1cm

$$\begin{matrix}
\bsmall
1
\underset{ {\aaa_{1}}}{\underparen{}}\infty
\underset{ {\aaa^*_{1}}}{\underparen{}}0
\underset{\underset{\circ}{\aaa_{2}}}{\underparen{}}\frac{1}{2}
\underset{\underset{\circ}{\aaa_{3}}}{\underparen{}}1

\esmall
&
1
\underset{ {\aaa_{1}}}{\underparen{}}\infty
\underset{ {\aaa^*_{1}}}{\underparen{}}0
\underset{ {\aaa_{2}}}{\underparen{}}\frac{1}{3}
\underset{ {\aaa_{3}}}{\underparen{}}\frac{1}{2}
\underset{ {\aaa^*_{3}}}{\underparen{}}\frac{2}{3}
\underset{ {\aaa^*_{2}}}{\underparen{}}1

&
1
\underset{ {\aaa_{1}}}{\underparen{}}\infty
\underset{ {\aaa^*_{1}}}{\underparen{}}0
\underset{ {\aaa_{2}}}{\underparen{}}\frac{1}{3}
\underset{ {\aaa^*_{2}}}{\underparen{}}\frac{2}{5}
\underset{ {\aaa_{3}}}{\underparen{}}\frac{1}{2}
\underset{ {\aaa^*_{3}}}{\underparen{}}1

\\
\begin{tikzpicture}[scale=3]
\scriptsize
\draw (0,1)--(0,0) node [very near start, right] {$\aaa_1^*$} node [below] {$0$}node (start) {} arc (180:0:0.25000000)  node [midway, above] {$\aaa_{2}$} node [below]{$\frac{1}{2}$} node (start) {} (2/5,1/5) node {$\circ$};
 \draw (start) node{}
node (start) {} arc (180:0:1/4)  node [midway, above] {$\aaa_{3}$} node (start) {} (3/5,1/5) node {$\circ$};
 \draw (start) node{} node [below] {$1$} -- (1,1) node [very near end, left] {$\aaa_1$};
\end{tikzpicture}
&
\begin{tikzpicture}
\scriptsize
\input fareydraw/cercle0\string_10.tex
\end{tikzpicture}
&
\begin{tikzpicture}
\scriptsize
\input fareydraw/cercle1\string_10.tex
\end{tikzpicture}
\\
\Gamma_0(5)&\Gamma_0(6)&\Gamma_0(6)
\end{matrix}
$$

$$
\Gamma_0(7) :
1
\underset{ {\aaa_{1}}}{\underparen{}}\infty
\underset{ {\aaa^*_{1}}}{\underparen{}}0
\underset{\underset{\bullet}{\aaa_{2}}}{\underparen{}}\frac{1}{2}
\underset{\underset{\bullet}{\aaa_{3}}}{\underparen{}}1

,\quad\Gamma_0(8) :
1
\underset{ {\aaa_{1}}}{\underparen{}}\infty
\underset{ {\aaa^*_{1}}}{\underparen{}}0
\underset{ {\aaa_{2}}}{\underparen{}}\frac{1}{4}
\underset{ {\aaa^*_{2}}}{\underparen{}}\frac{1}{3}
\underset{ {\aaa_{3}}}{\underparen{}}\frac{1}{2}
\underset{ {\aaa^*_{3}}}{\underparen{}}1

,\quad\Gamma_0(9) :
1
\underset{ {\aaa_{1}}}{\underparen{}}\infty
\underset{ {\aaa^*_{1}}}{\underparen{}}0
\underset{ {\aaa_{2}}}{\underparen{}}\frac{1}{3}
\underset{ {\aaa^*_{2}}}{\underparen{}}\frac{1}{2}
\underset{ {\aaa_{3}}}{\underparen{}}\frac{2}{3}
\underset{ {\aaa^*_{3}}}{\underparen{}}1

$$

\begin{equation*}
\Gamma_0(10) :
1
\underset{ {\aaa_{1}}}{\underparen{}}\infty
\underset{ {\aaa^*_{1}}}{\underparen{}}0
\underset{\underset{\circ}{\aaa_{2}}}{\underparen{}}\frac{1}{3}
\underset{ {\aaa_{3}}}{\underparen{}}\frac{2}{5}
\underset{ {\aaa_{4}}}{\underparen{}}\frac{1}{2}
\underset{ {\aaa^*_{4}}}{\underparen{}}\frac{3}{5}
\underset{ {\aaa^*_{3}}}{\underparen{}}\frac{2}{3}
\underset{\underset{\circ}{\aaa_{5}}}{\underparen{}}1

,\quad
1
\underset{ {\aaa_{1}}}{\underparen{}}\infty
\underset{ {\aaa^*_{1}}}{\underparen{}}0
\underset{\underset{\circ}{\aaa_{2}}}{\underparen{}}\frac{1}{3}
\underset{ {\aaa_{3}}}{\underparen{}}\frac{2}{5}
\underset{ {\aaa^*_{3}}}{\underparen{}}\frac{3}{7}
\underset{ {\aaa_{4}}}{\underparen{}}\frac{1}{2}
\underset{ {\aaa^*_{4}}}{\underparen{}}\frac{2}{3}
\underset{\underset{\circ}{\aaa_{5}}}{\underparen{}}1

\end{equation*}
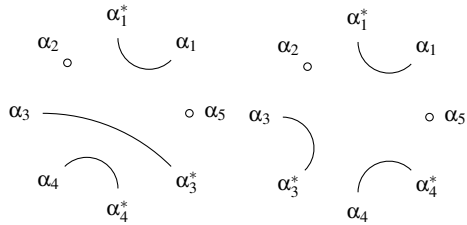
\begin{center}
\begin{tikzpicture}
\bsmall
\input fareydraw/cercle0\string_10.tex
\esmall
\end{tikzpicture}
\begin{tikzpicture}
\bsmall
\input fareydraw/cercle1\string_10.tex
\esmall
\end{tikzpicture}
\end{center}
\begin{equation*}
\Gamma_0(11) :
1
\underset{ {\aaa_{1}}}{\underparen{}}\infty
\underset{ {\aaa^*_{1}}}{\underparen{}}0
\underset{ {\aaa_{2}}}{\underparen{}}\frac{1}{3}
\underset{ {\aaa_{3}}}{\underparen{}}\frac{1}{2}
\underset{ {\aaa^*_{2}}}{\underparen{}}\frac{2}{3}
\underset{ {\aaa^*_{3}}}{\underparen{}}1

\end{equation*}
\begin{equation*}
\Gamma_0(12) :
1
\underset{ {\aaa_{1}}}{\underparen{}}\infty
\underset{ {\aaa^*_{1}}}{\underparen{}}0
\underset{ {\aaa_{2}}}{\underparen{}}\frac{1}{6}
\underset{ {\aaa^*_{2}}}{\underparen{}}\frac{1}{5}
\underset{ {\aaa_{3}}}{\underparen{}}\frac{1}{4}
\underset{ {\aaa_{4}}}{\underparen{}}\frac{1}{3}
\underset{ {\aaa_{5}}}{\underparen{}}\frac{1}{2}
\underset{ {\aaa^*_{5}}}{\underparen{}}\frac{2}{3}
\underset{ {\aaa^*_{4}}}{\underparen{}}\frac{3}{4}
\underset{ {\aaa^*_{3}}}{\underparen{}}1

,\quad
1
\underset{ {\aaa_{1}}}{\underparen{}}\infty
\underset{ {\aaa^*_{1}}}{\underparen{}}0
\underset{ {\aaa_{2}}}{\underparen{}}\frac{1}{6}
\underset{ {\aaa^*_{2}}}{\underparen{}}\frac{1}{5}
\underset{ {\aaa_{3}}}{\underparen{}}\frac{1}{4}
\underset{ {\aaa^*_{3}}}{\underparen{}}\frac{2}{7}
\underset{ {\aaa_{4}}}{\underparen{}}\frac{1}{3}
\underset{ {\aaa^*_{4}}}{\underparen{}}\frac{2}{5}
\underset{ {\aaa_{5}}}{\underparen{}}\frac{1}{2}
\underset{ {\aaa^*_{5}}}{\underparen{}}1

\end{equation*}
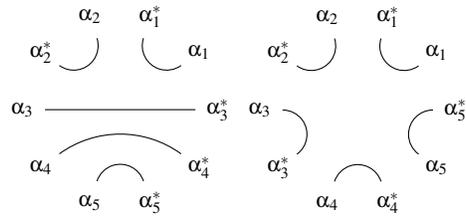
\begin{center}
\bsmall
\begin{tikzpicture}
\input fareydraw/cercle0\string_12.tex
\end{tikzpicture}
\begin{tikzpicture}
\input fareydraw/cercle1\string_12.tex
\end{tikzpicture}
\esmall
\end{center}
On voit sur ce dernier \farey que $X_0(12)$ est de genre 0, n'a pas de points elliptiques et
a 4 pointes (les pointes $p$ autres que $0$ sont en bijection
avec les chemins $\aaa \aaa^*$ allant de $0$ à $0$ en passant par $p$).
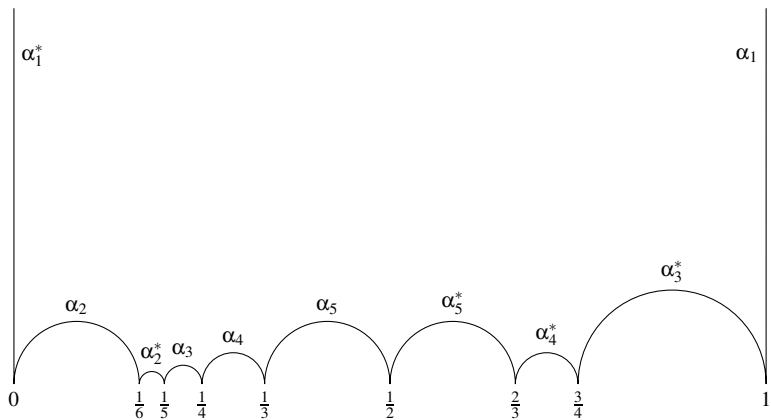
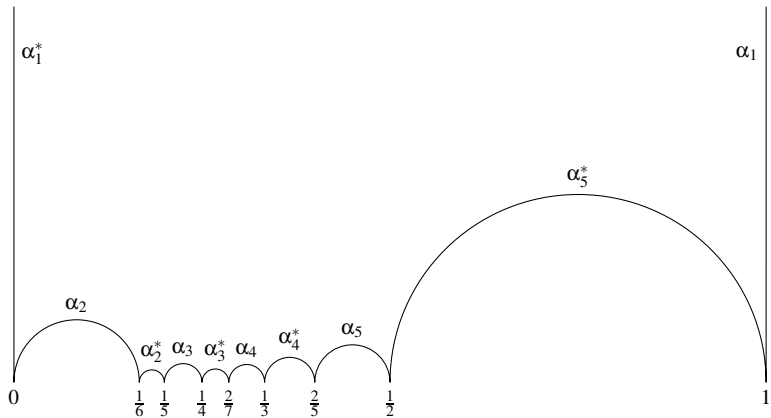
\begin{center}
\input polygondraw/0\string_12.tex
\end{center}
\begin{center}
\input polygondraw/2\string_12.tex
\end{center}

\begin{equation*}
\Gamma_0(13) :
1
\underset{ {\aaa_{1}}}{\underparen{}}\infty
\underset{ {\aaa^*_{1}}}{\underparen{}}0
\underset{\underset{\bullet}{\aaa_{2}}}{\underparen{}}\frac{1}{3}
\underset{\underset{\circ}{\aaa_{3}}}{\underparen{}}\frac{1}{2}
\underset{\underset{\circ}{\aaa_{4}}}{\underparen{}}\frac{2}{3}
\underset{\underset{\bullet}{\aaa_{5}}}{\underparen{}}1

\end{equation*}
\begin{center}
\input polygondraw/2\string_13.tex
\end{center}
\begin{equation*}
\Gamma_0(14) :
1
\underset{ {\aaa_{1}}}{\underparen{}}\infty
\underset{ {\aaa^*_{1}}}{\underparen{}}0
\underset{ {\aaa_{2}}}{\underparen{}}\frac{1}{4}
\underset{ {\aaa_{3}}}{\underparen{}}\frac{2}{7}
\underset{ {\aaa_{4}}}{\underparen{}}\frac{1}{3}
\underset{ {\aaa_{5}}}{\underparen{}}\frac{2}{5}
\underset{ {\aaa^*_{4}}}{\underparen{}}\frac{3}{7}
\underset{ {\aaa^*_{3}}}{\underparen{}}\frac{1}{2}
\underset{ {\aaa^*_{2}}}{\underparen{}}\frac{2}{3}
\underset{ {\aaa^*_{5}}}{\underparen{}}1

,\quad
1
\underset{ {\aaa_{1}}}{\underparen{}}\infty
\underset{ {\aaa^*_{1}}}{\underparen{}}0
\underset{ {\aaa_{2}}}{\underparen{}}\frac{1}{4}
\underset{ {\aaa^*_{2}}}{\underparen{}}\frac{7}{27}
\underset{ {\aaa_{3}}}{\underparen{}}\frac{2}{7}
\underset{ {\aaa^*_{3}}}{\underparen{}}\frac{13}{43}
\underset{ {\aaa_{4}}}{\underparen{}}\frac{10}{33}
\underset{ {\aaa_{5}}}{\underparen{}}\frac{1}{3}
\underset{ {\aaa^*_{4}}}{\underparen{}}\frac{2}{5}
\underset{ {\aaa^*_{5}}}{\underparen{}}1

\end{equation*}
\begin{center}
\begin{tikzpicture}
\scriptsize
\input fareydraw/cercle0\string_14.tex
\end{tikzpicture}
\begin{tikzpicture}
\bsmall
\input fareydraw/cercle1\string_14.tex
\esmall
\end{tikzpicture}
\end{center}
Les arcs $\aaa_4$ and $\aaa_5$ du deuxième \farey se projettent dans $X_0(14)$ en
des lacets et donnent une base symplectique de $H_1(X_0(14),\ZZ)$ pour le
produit d'intersection. Les arcs $\aaa_4$ et $\aaa_5$ sont de type hyperbolique,
les arcs $\aaa_1$, $\aaa_2$ et $\aaa_3$ sont de type parabolique.

\foreach \n in {15,16,...,29,31,...,33,35,37}{

\begin{table}
\caption{$N=\n$}
\begin{center}
\begin{tabular}{m{3.5cm}m{10cm}}
\begin{tikzpicture}[scale=1]
\bsmall
\input fareydraw/cercle0\string_\n.tex
\esmall
\end{tikzpicture}
&
$
\input fareydraw/0\string_\n.tex
$
\end{tabular}
\end{center}
\begin{center}
\begin{tabular}{m{3.5cm}m{10cm}}
\begin{tikzpicture}[scale=1]
\bsmall
\input fareydraw/cercle1\string_\n.tex
\esmall
\end{tikzpicture}
&
\ifnum\n=24\bsmall\fi
\ifnum\n=28\bsmall\fi
\ifnum\n>34\bsmall\fi
$
\input fareydraw/1\string_\n.tex
$
\ifnum\n=24\esmall\fi
\ifnum\n=28\esmall\fi
\ifnum\n>34\esmall\fi
\end{tabular}
\end{center}
\end{table}
}
%%si moins de dessins
%%\endinput
\clearpage
Nous ne montrons désormais que la représentation des involutions.
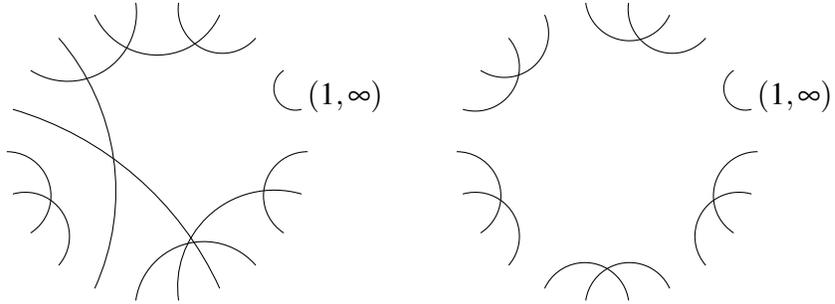
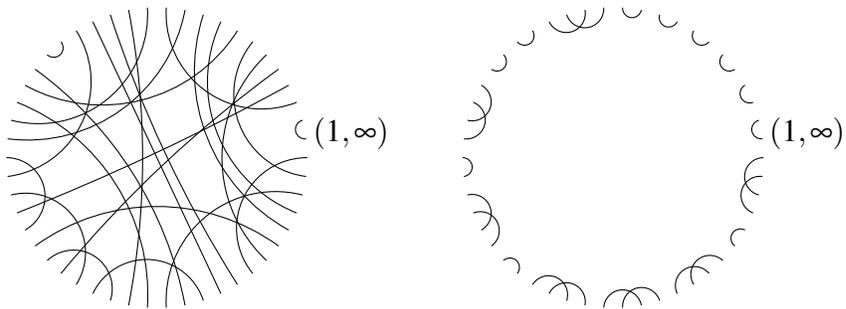
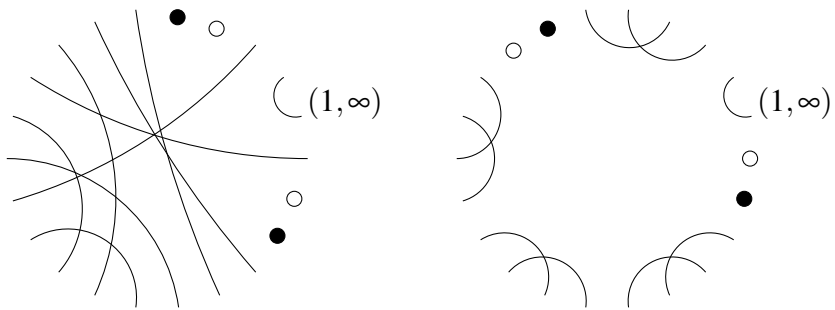
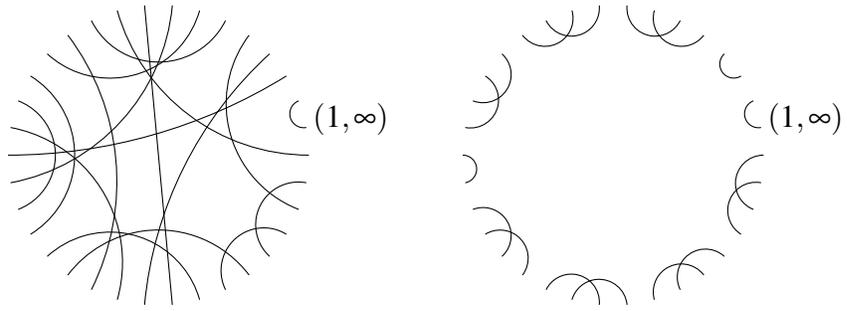
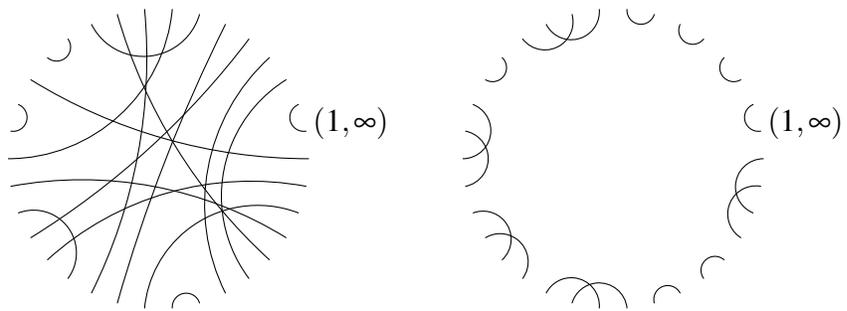
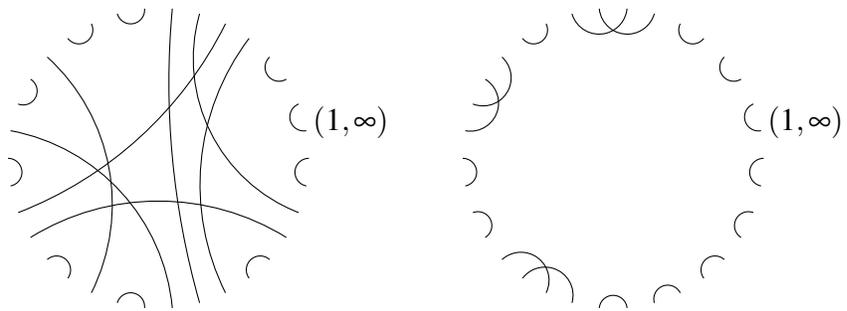
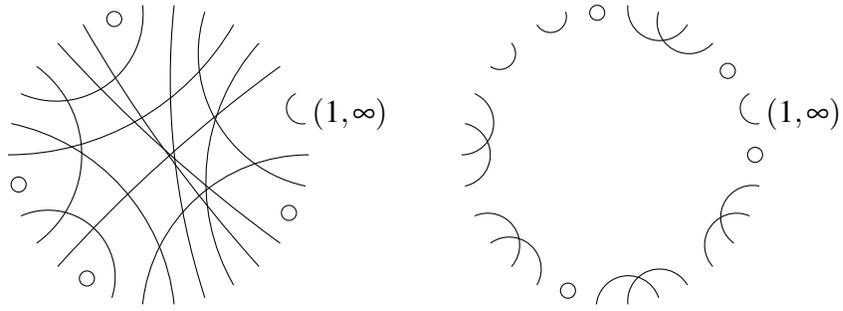
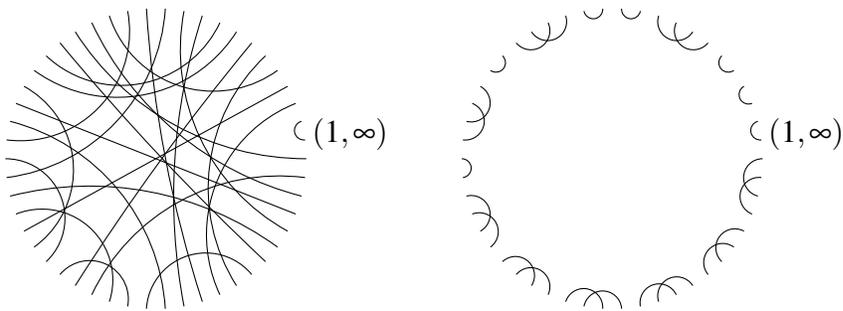
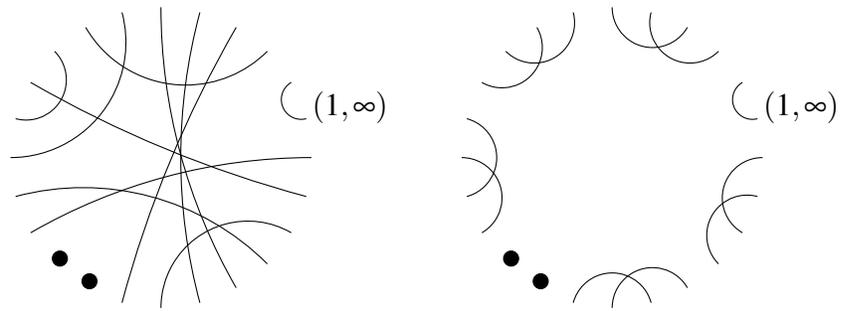
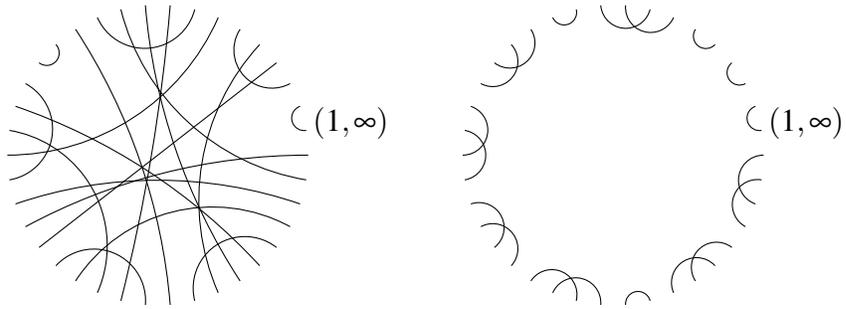
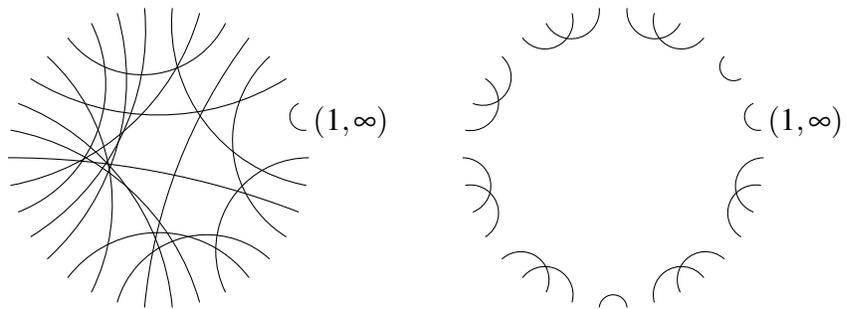
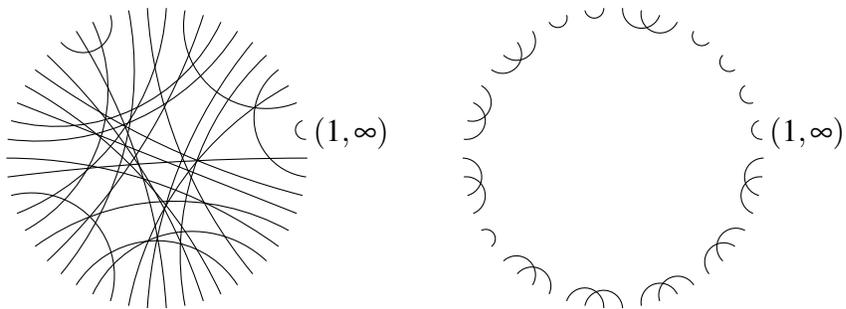
\foreach \n in {59,..., 70}{
$$N=\n$$
\begin{center}
\begin{tikzpicture}[scale=2]
\input fareydraw/cercle0\string_\n.tex
\end{tikzpicture}
\begin{tikzpicture}[scale=2]
\input fareydraw/cercle1\string_\n.tex
\end{tikzpicture}
\end{center}
\hrule
}

 % !TEX encoding = IsoLatin


\begin{thebibliography}{XX}

\bibitem{bost} Jean-Benoit Bost,
Introduction to Compact Riemann Surfaces, Jacobians, and Abelian Varieties
in From number theory to physics (Les Houches, 1989), Springer, pp 64-211.

\bibitem{congruence}
Ann Dooms, Eric Jespers, Alexander Konovalov and Helena Verrill.
\kbd{Congruence} - Congruence subgroups of $\Sl_2(\ZZ)$,
Version 1.2.2; 2018 (\url{https://www.gap-system.org/Packages/congruence.html}).

\bibitem{cremona}
John Cremona, Algorithms for modular elliptic curves (2nd ed.), Cambridge
University Press, Cambridge, 1997.

\bibitem{kulkarni}
Ravi S. Kulkarni, An Arithmetic-Geometric Method in
the Study of the Subgroups of the Modular Group,
American Journal of Mathematics, Vol. 113, No. 6. (Dec., 1991), pp. 1053-1133.

\bibitem{KurtLong08}
Chris A. Kurth et Ling Long, Computations with finite index subgroups of
$\Psl_2(\ZZ)$ using Farey symbols. Advances in algebra and combinatorics,
225-242, World Sci. Publ., Hackensack, NJ, 2008.

\bibitem{maskit} Bernard Maskit, On Poincare's Theorem for Fundamental Polygons,
Advances in Math., 7 (1971). 219-230.

\bibitem{pari} The PARI~Group, PARI/GP version \texttt{2.11.1}, Univ.
Bordeaux, 2018, \url{http://pari.math.u-bordeaux.fr}.

\bibitem{poincare} Henri Poincaré, Théorie des groupes fuchsiens, Acta Math.
Volume 1 (1882), 1-62.

\bibitem{PS} Robert Pollack, Glenn Stevens,
Overconvergent modular symbols and p-adic L-functions,
Annales scientifiques de l'ENS 44, fascicule 1 (2011), 1-42.

%\bibitem{sage}
%SageMath, the Sage Mathematics Software System (Version 8.6),
%   The Sage Developers, 2018, \url{https://www.sagemath.org}.

\bibitem{siegel}
Carl L. Siegel, Topics in Complex Function Theory,
vol I, Elliptic Functions and Uniformization Theory, Interscience Tracts in
Pure and Applied Mathematics, No. 25, Vol. I, 1969.

\bibitem{verrill}
Helena Verrill, Fundamental domain drawer, 2000,
archived at \url{https://wstein.org/Tables/fundomain/}.
\end{thebibliography}
\end{document}